\documentclass[10pt,a4paper,twoside]{amsart}

\usepackage[utf8]{inputenc}
\usepackage[english]{babel}
\usepackage[T1]{fontenc}
\usepackage{microtype, fancyhdr}
\usepackage{listings}
\usepackage[a4paper, hmarginratio=1:1]{geometry}  
\usepackage{kerkis}
\usepackage{charter}

\usepackage{multicol}
\usepackage{amsfonts, amsmath, amsthm, amssymb, mathrsfs, amscd, mathtools}
\usepackage[all]{xy}
\numberwithin{equation}{section}	
\usepackage{faktor}
\allowdisplaybreaks

\usepackage{url, enumitem}
\usepackage[noadjust]{cite}

\usepackage{graphicx}
\graphicspath{{Figures/}}

\usepackage{hyperref}

\theoremstyle{plain}
\newtheorem{theorem}{Theorem}[section]
\newtheorem{proposition}[theorem]{Proposition}
\newtheorem{corollary}[theorem]{Corollary}
\newtheorem{lemma}[theorem]{Lemma}

\theoremstyle{remark}
\newtheorem{remark}[theorem]{Remark}

\theoremstyle{definition}
\newtheorem{definition}[theorem]{Definition}

\makeatletter
\def\@map#1#2[#3]{\mbox{$#1 \colon\thinspace #2 \longrightarrow #3$}}
\def\map#1#2{\@ifnextchar [{\@map{#1}{#2}}{\@map{#1}{#2}[#2]}}
\g@addto@macro\@floatboxreset\centering
\makeatother

\newcommand\Nn{\mathbb{N}}
\newcommand\Zz{\mathbb{Z}}

\newcommand\Rr{\mathbb{R}}

\DeclareMathOperator{\Aut}{Aut}

\newcommand\BB[1]{B_{#1}}	
\newcommand\BC[2]{B_{#1}[#2]}	

\newcommand{\aut}[1]{\ensuremath{\operatorname{\text{Aut}}\left({#1}\right)}}
\newcommand{\aff}[1]{\ensuremath{\operatorname{\text{Aff}}\left({#1}\right)}}

\newcommand\sig[1]{\sigma_{\hspace{-0.3ex}#1}^{\null}} 
\newcommand\sigg[2]{\sigma_{\hspace{-0.3ex}#1}^{#2}}	

\newcommand\Sp[2]{\mathrm{Sp}_{#1}(#2)}	
\newcommand\PSp[2]{\mathrm{PSp}_{#1}(#2)} 

\newcommand\ii{i} 
\newcommand\nn{n} 
\newcommand\nno{n-1} 
\newcommand\mm{m} 

\usepackage{listings}
\lstdefinelanguage{GAP}{%
  morekeywords={%
    Assert,Info,IsBound,QUIT,%
    TryNextMethod,Unbind,and,break,%
    continue,do,elif,%
    else,end,false,fi,for,%
    function,if,in,local,%
    mod,not,od,or,%
    quit,rec,repeat,return,%
    then,true,until,while%
  },%
  sensitive,%
  morecomment=[l]\#,%
  morestring=[b]",%
  morestring=[b]',%
}[keywords,comments,strings]

\usepackage[T1]{fontenc}
\usepackage[variablett]{lmodern}
\usepackage{xcolor}
\begin{document}
\title{Congruence subgroups of braid groups and crystallographic quotients. Part II}

\author[Bellingeri]{Paolo Bellingeri}
\address{Normandie Univ, UNICAEN, CNRS, LMNO, 14000 Caen, France}
\email{paolo.bellingeri@unicaen.fr}

\author[Damiani]{Celeste Damiani}
\address{Fondazione Istituto Italiano di Tecnologia, Genova, Italy}
\email{celeste.damiani@iit.it}

\author[Ocampo]{Oscar Ocampo}
\address{Universidade Federal da Bahia, Departamento de Matem\'atica - IME, CEP:~40170-110 - Salvador, Brazil}
\email{oscaro@ufba.br}

\author[Stylianakis]{Charalampos Stylianakis}
\address{University of the Aegean, Department of mathematics, Karlovasi, 83200, Samos, Greece}
\email{stylianakisy2009@gmail.com}



\date{\today}

\begin{abstract}
Following previous work on congruence subgroups and crystallographic braid groups, we study the lower central series of congruence braid groups related to the braid group $B_3$, showing in particular that the corresponding quotients are almost crystallographic.
\end{abstract}

\maketitle

\section{Introduction}

In this paper we continue the study, initiated in~\cite{BDOS:2024}, of the relationship between congruence subgroups of the braid group $B_n$ and crystallographic braid groups respectively introduced Arnol'd~\cite{Arnold:1968} and Tits~\cite{JTits:1966}.

As recalled in ~\cite{BDOS:2024} both families can be defined in  a general settings but recently they have been intensively studied in the specific case of braid groups and relatives,  see for instance~\cite{Brendle-Margalit:2018, Stylianakis:2018, ABGH, Nakamura:2021, KordekMargalit:2022, bloomquist2023quotients,Banerjee-Huxford}  for congruence subgroups of braid groups and~\cite{ Goncalves-Guaschi-Ocampo:2017, BeckMarin:2020, Goncalves-Guaschi-Ocampo-Pereiro, Bellingeri-Guaschi-Makri, Cerqueira-Ocampo} for crystallographic braid groups. 

In the context of groups of matrices, a \emph{congruence subgroup} of a matrix group with integer entries 
is a subgroup defined as the kernel of the mod $m$ reduction of a linear group.
The notion of congruence subgroups can be generalised for mapping class groups using the action on the homology of the surface. The mapping classes acting trivially on the homology forms the Torelli group of the surface while congruence subgroups can be obtained  considering the homology of the surface mod $p$.
In the case of the  braid group $\BB\nn$  a similar (symplectic) representation is obtained 
considering integral Burau representation, whose kernel gives the analogous of Torelli group for braid  \cite[Proposition 2.1]{Gambaudo-Ghys:2005}. The level $m$ congruence subgroup of  $B_n$ is therefore defined as the kernel 
of the Burau representation on $\Zz/m \Zz$ coefficients and denoted by $\BC\nn\mm$.

On the other hand crystallographic groups are extensions of finite groups by an abelian free group of finite rank. Gon\c{c}alves-Guaschi-Ocampo \cite{Goncalves-Guaschi-Ocampo:2017} showed that the group $\faktor{\BB\nn}{[P_n,P_n]}$ is crystallographic, where $P_n$ is the pure braid subgroup of $P_n$ and where $[G,G]$ denotes the commutator of $G$. Motivated by the fact that   $\BC\nn2$ coincides with $P_n$
(see \cite{Arnold:1968} or later on the paper) in~\cite{BDOS:2024} we established some (iso)-morphisms between crystallographic braid groups
and the corresponding quotients of congruence 
braid groups.  Actually, Gon\c{c}alves-Guaschi-Ocampo proved that for $n,k\geq 3$ then $\faktor{\BB\nn}{\Gamma_k(\BC{\nn}{2})}$ is an almost-crystallographic group \cite{LIMAGONCALVES2019160}, where $\Gamma_k(\BC{\nn}{2})$ is the $k$-th subgroup of the lower central series of $\BC\nn{2}$ (see next section for the definition of almost-crystallographic group and of the lower central series of a group).

In this paper, we consider $p$ to be an odd prime or $p=4$, and we focus on the whole central series of $\Gamma_k(\BC{n}{p})$ in the specific case  $n=3$. 
The three strands braid group $\BB3$, can both be used as a toy model to understand properties of braid groups on more strands, and as a special case among braid groups, as it presents a number of properties that single it out from the rest of the family.
Several group theoretic and combinatorial properties can be proved leveraging Artin's presentation for~$\BB3$, easier to manipulate than the presentation for the general braid group~$\BB\nn$. Besides, this group appears as the fundamental group of the trefoil knot
and mapping class group of a surface of genus one
with one boundary component. All these facts participate to the richness of the representation theory of $B_3$~\cite{Tuba-Wenzl:2000}, which is also motivated by fascinating physical and topological applications~\cite{Bruillard-Plavnik-Rowell:2015}.
However we think that our results could open possible leads to the generic case.

We have the following.

\begin{theorem}
Let $p$ be an odd prime or $p=4$. Then:
\begin{enumerate}
\label{main_t}
    \item The group $\faktor{\BB3}{[\BC{3}{p},\BC{3}{p}]}$ is crystallographic.
    \item For $k\geq 2$ the group $\faktor{\BB3}{\Gamma_k(\BC{3}{p})}$  is almost-crystallographic.
\end{enumerate}
\end{theorem}

The first part of Theorem~\ref{main_t} is a consequence of Corollary~\ref{b3p} and Theorem~\ref{T:34} if $p$ is prime; it is proved in Proposition~\ref{P:b34cryst} when~$p=4$. The second part of Theorem \ref{main_t} is proven in Theorem~\ref{T:alm_cryst}.

Recall  that a torsion free crystallographic group is called a Bieberbach group (see also Definition 2.1). 
By specializing to the case $p=3$, the following is a consequence of Theorem~\ref{torsionfree_b3}.

\begin{theorem}
    The group $\faktor{\BB3}{[\BC3{3},\BC3{3}]}$ is a Bieberbach group with holonomy group $\PSp2{\Zz/3\Zz}$.
    \label{b3p_Bieberbach}
\end{theorem}

This paper is organized as follows. In Section~\ref{S:prelim}, we briefly recall the basic definitions used throughout the paper. In particular,  we review the symplectic representation and certain distinguished subgroups of the braid group, including the center, the pure braid group, and the congruence braid subgroups. We also recall the relevant notions and results concerning crystallographic and almost-crystallographic groups. Section~\ref{S:presentation} is devoted to the proof of the first part of Theorem~\ref{main_t} as well as to Theorem~\ref{b3p_Bieberbach}. 
A key ingredient 
is a presentation of $\BC3{p}$ where $p$ is prime or $p=4$
which yields a direct product decomposition $\BC3{p}\cong \Zz \times F$, with $F$ a  free group. 
This decomposition is established in Proposition~\ref{P:b34} and Theorem~\ref{b3p_presentation}. 
In Section~\ref{S:almost}, we establish the relationship between almost-crystallographic structures and congruence subgroups, which completes the proof of the second part of Theorem~\ref{main_t}. 
Finally, in Appendix~\ref{S:Appendix} 
we examine the algebraic structure of the finite group $\rho_4(B_3)$. This group is of independent interest, and also plays a key role in proving that the quotient $\faktor{B_3}{[B_3[4], B_3[4]]}$ is a crystallographic group of dimension~6 with holonomy group~$S_4$.


\subsection*{Acknowledgments}
P.~B. was partially supported by the ANR project AlMaRe (ANR-19-CE40-0001). C.~D.  is a member of GNSAGA of INdAM. O.~O. acknowledges the hospitality of the \textit{Laboratoire de Math\'ematiques Nicolas Oresme} (Universit\'e de Caen Normandie) from August 2023 to January 2024, where part of this project was developed, and the partial support from CAPES (Programa CAPES-Print, Process no. 88887.835402/2023-00) and from the National Council for Scientific and Technological Development (CNPq) through a \textit{Bolsa de Produtividade} (305422/2022-7). O.~O. wishes to express his gratitude to Eliane Santos, all HCA staff, Bruno Noronha, Luciano Macedo, M\'arcio Isabella, Andreia de Oliveira Rocha, Andreia Gracielle Santana, and Ednice de Souza Santos for their valuable help since July 2024. C.~S. would like to thank Peter Patzt for useful conversations on the structure of congruence subgroups of braid groups. The authors are grateful to the anonymous referee for a careful reading and valuable suggestions that improved the manuscript.

\section{Preliminaries}\label{S:prelim}

Let us begin by recalling a few notations and definitions.

\subsection{Center, lower central and derived series}\label{S:defs}

Let $G$ be a group. We denote its center by $Z(G)$. We recall that its \emph{lower central series} $\{\Gamma_{k}(G)\}_{k\in \mathbb{N}}$ is defined by $\Gamma_{1}(G)=G$, and $\Gamma_{k}(G)=[\Gamma_{k-1}(G),G]$ for all $k\geq 2$ (if $H$ and $K$ are subgroups of $G$, $[H,K]$ is defined to be the subgroup of $G$ normally generated by the commutators of the form $[h,k]=hkh^{-1}k^{-1}$, where $h\in H$ and $k\in K$). The \emph{derived series} of $G$, $\{ G^{(i)} \}_{i\in \Nn\cup \{0\}}$, is defined inductively by $G^{(0)} = G$, and $G^{(i)} =
[G^{(i-1)},\, G^{(i-1)}]$ for all $i \in \Nn$.

It is worth noting that $\Gamma_{2}(G)=G^{(1)}$ is the commutator subgroup of $G$, and the groups $\Gamma_{k}(G)$, $G^{(k)}$ are normal subgroups of $G$, for all $k\in \mathbb{N}$.

\subsection{Braid group and congruence subgroups}\label{Ss:cong}

In the following the  \emph{braid group} $\BB\nn$ is defined as the mapping class group of the $n$-punctured disk. We recall that it admits the following presentation~\cite{Artin:1925}:
\[
\bigg\langle \sig1, \ldots\, , \sig\nno \ \bigg\vert \ 
\begin{matrix}
\sig{i} \sig j = \sig j \sig{i} 
&\text{for} &\vert  i-j\vert > 1\\
\sig{i} \sig j \sig{i} = \sig j \sig{i} \sig j 
&\text{for} &\vert  i-j\vert  = 1
\end{matrix}
\ \bigg\rangle.
\]
Recall that the \emph{full twist} $\Delta_n^2 = (\sig1 \cdots \sigma_{n-1})^n$ generates the center $Z(\BB\nn)$.

The reducible Burau representation is a linear representation $b_{t}\colon \BB\nn \to \mathrm{GL}_n(\mathbb{Z}[t^{\pm 1}])$ defined as follows
\[
b_{t}(\sig\ii) = I_{i-1} \oplus \begin{pmatrix}
1-t&t\\	
1&0
\end{pmatrix} \oplus I_{n-i-1},
\]
where $I_k$ is the identity matrix of dimension $k$. By evaluating it in $t=-1$  we obtain the integral Burau representation $b_{-1}\colon \BB\nn \to \mathrm{GL}_n(\mathbb{Z})$. By the work of Gambaudo-Ghys, the image of $b_{-1}$ is a subgroup of a symplectic group~\cite[Proposition2.1]{Gambaudo-Ghys:2005}. Consequently, we obtain
\begin{equation*}
 \rho\colon \BB\nn \rightarrow 
\begin{cases} 
\Sp{\nn-1}{\Zz} \mbox{ for } \nn \mbox{ odd}, \\
(\Sp{\nn}{\Zz})_u \mbox{ for } \nn \mbox{ even}, \\
\end{cases}
\end{equation*}
where $(\Sp{\nn}{\Zz})_u$ is a subgroup of $\Sp{\nn}{\Zz}$ fixing one vector. The reduction map $\Zz \to \Zz/m\Zz$ induces a map
\begin{equation}
\label{eqn:rhom}
 \rho_m\colon \BB\nn \rightarrow 
\begin{cases} 
\Sp{\nn-1}{\Zz/m\Zz} \mbox{ for } \nn \mbox{ odd}, \\
(\Sp{\nn}{\Zz/m\Zz})_u \mbox{ for } \nn \mbox{ even}. \\
\end{cases}
\end{equation}
The kernel of $\rho_m$ is denoted by $\BC\nn{m}$ and it is called \emph{level $m$ congruence subgroup} of $\BB\nn$. When $n=3$ the group $\rho_3(B_3)$ is isomorphic to $\Sp{2}{\Zz/3\Zz}$ \cite[Theorem~1]{ACampo:1979} and the corresponding matrices of $\sig1, \sig2$ are as follows \cite[Section 2.2]{bloomquist2023quotients}:

$$
\sig1: 
\begin{pmatrix}
1&1\\	
0&1
\end{pmatrix}
\quad
\textrm{ and }
\quad
\sig2: 
\begin{pmatrix}
1&0\\	
-1&1
\end{pmatrix}.
$$

The whole construction above is detailed in \cite[Section 2.2]{BDOS:2024}.

\subsection{Pure braid groups} 
\label{Ss:pure}

Recall that the pure braid group $P_n$
is the kernel of the map from $B_n$ to $S_n$ associating to any braid the corresponding permutation.

By definition, we have that $\sigg{i}{2}$ is an element of~$P_n$. For $i < j$, let us define braid words
\[ A_{i,j} = (\sigma_{j-1} \sigma_{j-2} ... \sigma_{i+1}) \sigg{i}{2} (\sigma_{j-1} \sigma_{j-2} ... \sigma_{i+1})^{-1}. \]
The group $P_n$ has a presentation given by generators $A_{i,j}$ with the following defining relations, see for instance \cite[Section 1.2]{BirmanBrendle:2005}:
\[
\begin{matrix*}[l]
A^{-1}_{r,s} A_{i,j} A_{r,s}  = A_{i,j} & 1 \leq r<s<i<j \leq n \: \mathrm{or} \: 1 \leq i<r<s<j \leq n, \\
A^{-1}_{r,s} A_{i,j} A_{r,s}  = A_{r,j} A_{i,j} A^{-1}_{r,j}  & 1 \leq r<s=i<j \leq n, \\
A^{-1}_{r,s} A_{i,j} A_{r,s}  = (A_{i,j} A_{s,j})A_{i,j}(A_{i,j} A_{s,j})^{-1} & 1 \leq r=i<s<j \leq n, \\
A^{-1}_{r,s} A_{i,j} A_{r,s}  = [(A_{r,j} A_{s,j} A^{-1}_{r,j} A^{-1}_{s,j}),A_{i,j}]A_{i,j} & 1 \leq r<i<s<j \leq n.
\end{matrix*}. 
\]

We can see that the full twist $\Delta_n^2$ is an element of $P_n$, since $\sigma_1 \sigma_2 \cdots \sigma_{n-1}$ represents the $n$-cycle $(1,2,\cdots, n)$ and hence $\Delta_n^2 = (\sigma_1 \sigma_2 \cdots \sigma_{n-1})^n$ has associated permutation $(1,2,\cdots, n)^n$ which is trivial in $S_n$.

Moreover, the pure braid group $P_n$ is isomorphic to the congruence subgroup $B_n[2]$ and $\mathcal{N}_n(\sigma_1^2)$, the normal closure of $\sigma_1^2$ in $B_n$. For $m\geq 3$, the index of $\mathcal{N}_n(\sigg{1}{m})$, the normal closure of $\sigma_1^m$ in $B_n$, was determined by Coxeter in \cite{Coxeter:1957}, so in the following we will call these normal closures as \emph{Coxeter braid subgroups}. The general relation between congruence and Coxeter braid subgroups will be developed in a further paper. However, we note here that, according to~\cite[Lemma~2.1]{BDOS:2024} we have that for any $m\ge 2$, $\mathcal{N}_n(\sigg{1}{m})$ is a normal subgroup of $\BC\nn\mm$.

\subsection{Crystallographic groups and almost crystallographic groups}

In this subsection, we recall basic definitions and useful results about crystallographic groups and almost crystallographic groups. 

\subsubsection{Crystallographic groups}

We start by recalling some definitions and facts about crystallographic and Bieberbach groups, and the characterisation of crystallographic groups in terms of a representation that arises in certain group extensions whose kernel is a free Abelian group of finite rank and whose quotient is finite. We also recall some results concerning Bieberbach groups and the fundamental groups of flat Riemannian manifolds. For more details, see~\cite[Section~I.1.1]{Charlap},~\cite[Section~2.1]{Dekimpe} or~\cite[Chapter~3]{Wolf}.

Let $G$ be a Hausdorff topological group. A subgroup $H$ of $G$ is said to be \emph{discrete} if it is a discrete subset. If $H$ is a closed subgroup of $G$ then the quotient space $G/H$ admits the quotient 
topology for the canonical projection $\map{\pi}{G}[G/H]$, and we say that $H$ is \emph{uniform} if $G/H$ is compact.
From now on, we identify $\operatorname{\text{Aut}}(\mathbb{Z}^m)$ with $\operatorname{\text{GL}}(m,\mathbb{Z})$. 

\begin{definition}[{\cite[Sec.~2.1, p.~13]{Dekimpe}}]
A discrete, uniform subgroup $\Pi$ of $\Rr^m\rtimes \operatorname{\text{O}}(m,\Rr)\subseteq \operatorname{\text{Aff}}(\Rr^m)$ is said to be a \textit{crystallographic group} of dimension $m$. If in addition $\Pi$ is torsion free then $\Pi$ is called a \textit{Bieberbach group} of dimension $m$.
\end{definition}

If $\Phi$ is a group, an \emph{integral representation of rank $m$ of $\Phi$} is defined to be a homomorphism $\map{\Theta}{\Phi}[\operatorname{\text{Aut}}(\mathbb{Z}^m)]$. Two such representations are said to be \emph{equivalent} if their images are conjugate in $\operatorname{\text{Aut}}(\mathbb{Z}^m)$. We say that $\Theta$ is a \emph{faithful representation} if it is injective. We recall the following useful characterisation of crystallographic groups.

\begin{lemma}[{\cite[Lemma~8]{Goncalves-Guaschi-Ocampo:2017}} ]\label{L:cryst}
Let $\Pi$ be a group. Then $\Pi$ is a crystallographic group if and only if there exists an integer $m\in \mathbb N$, a finite group $\Phi$ and a short exact sequence of the form:
\begin{equation}\label{eq:SeqCrist2}
0\to \mathbb{Z}^m \to \Pi \stackrel{\zeta}{\longrightarrow} \Phi \to 1,
\end{equation}
such that the integral representation $\map{\Theta}{\Phi}[\operatorname{\text{Aut}}(\mathbb{Z}^m)]$ induced by conjugation on $\mathbb{Z}^m$ and defined by $\Theta(\phi)(x)=\pi x \pi^{-1}$ for all $x\in \mathbb{Z}^{m}$ and $\phi\in \Phi$, where $\pi\in \Pi$ is such that $\zeta(\pi)=\phi$, is faithful. 
\end{lemma}

If $\Pi$ is a crystallographic group, the integer $m$, the finite group $\Phi$ and the integral representation $\map{\Theta}{\Phi}[\operatorname{\text{Aut}}(\mathbb{Z}^{m})]$ appearing in the statement of Lemma~\ref{L:cryst} are called the \emph{dimension}, the \emph{holonomy group} and the \emph{holonomy representation} of  $\Pi$ respectively.

We recall the following two useful general results that will be applied to the study of crystallographic structures on quotients of the braid group by commutator subgroups of congruence subgroups.

\begin{theorem}[{\cite[Theorem~3.3]{BDOS:2024}}]
    Let $\phi\colon G\to F$ be a surjective homomorphism with $F$ a finite group. 
    Let $K$ denote the kernel of $\phi$. Suppose that there is a non-trivial element of the center of $G$ that does not belong to $K$.
    Then the representation $\eta\colon F\to \mathrm{Aut}\left(\faktor{K}{[K,K]}\right)$, induced from the action by conjugation of $\faktor{G}{[K,K]}$ on $\faktor{K}{[K,K]}$,  is not injective.
    \label{T:33}
\end{theorem}

In \cite{BDOS:2024} we considered the case where the holonomy representation defined in Lemma~\ref{L:cryst} is not injective and give conditions for the middle group to be a crystallographic group. 

\begin{theorem}[{\cite[Theorem~3.4]{BDOS:2024}}]
Consider the short exact sequence
$
1 \longrightarrow K \longrightarrow G \stackrel{p}{\longrightarrow} Q \longrightarrow 1
$
where $K$ is a free abelian group of finite rank and $Q$ is a finite group such that the representation $\varphi\colon Q \to Aut(K)$, induced from the action by conjugation, is not injective. 
Suppose that the group $p^{-1}(Ker(\varphi))$ is torsion free. 
Then $G$ is a crystallographic group with holonomy group $\faktor{Q}{Ker(\varphi)}$.
\label{T:34}
\end{theorem}

We now recall the connection between Bieberbach groups and manifolds, which motivates the study of torsion-free crystallographic groups. A Riemannian manifold $M$ is called \emph{flat} if it has zero curvature at every point. By the first Bieberbach Theorem, there is a correspondence between Bieberbach groups and fundamental groups of flat Riemannian manifolds without boundary (see~\cite[Theorem~2.1.1]{Dekimpe} and the paragraph that follows it). By~\cite[Corollary~3.4.6]{Wolf}, the holonomy group of a flat manifold $M$ is isomorphic to the group $\Phi$. In 1957, Auslander and Kuranishi proved that every finite group is the holonomy group of some flat manifold~(see~{\cite[Theorem~3.4.8]{Wolf} and~\cite[Theorem~III.5.2]{Charlap}}). It is well known that a flat manifold determined by a Bieberbach group $\Pi$ is orientable if and only if the integral representation $\map{\Theta}{\Phi}[\operatorname{\text{GL}}(m,\mathbb{Z})]$ satisfies $\operatorname{\text{Im}}(\Theta) \subseteq \operatorname{\text{SL}}(m,\mathbb{Z})$~\cite[Theorem~6.4.6 and Remark~6.4.7]{Dekimpe}.  This being the case, $\Pi$ is said to be an \emph{orientable Bieberbach group}.  

\subsubsection{Almost-crystallographic groups}

We recall briefly the definitions of almost-crystallographic and almost-Bieberbach groups, which are natural generalisations of crystallographic and Bieberbach groups, as well as a characterisation of almost-crystallographic groups. For more details about (almost-)crystallographic groups, see~\cite[Section~2.1]{Dekimpe}. 

Given a connected and simply-connected nilpotent Lie group $N$, the group $\aff{N}$ of affine transformations of $N$ is defined by $\aff{N}=N\rtimes \aut{N}$, and acts on $N$ by:
\begin{equation*}
\text{$(n,\phi)\cdot m = n\phi(m)$ for all $m,n\in N$ and $\phi\in \aut{N}$}.
\end{equation*}

\begin{definition}[{\cite[Sec.~2.2, p.~15]{Dekimpe}}]\label{D:almostc}
Let $N$ be a connected, simply-connected nilpotent Lie group, and consider a maximal compact subgroup $C$ of $\aut{N}$. A uniform discrete subgroup $E$ of $N\rtimes C$ is called an \emph{almost-crystallographic group}, and its dimension is defined to be that of $N$. A torsion-free, almost-crystallographic group is called an \emph{almost-Bieberbach group}, and the quotient space $E\backslash N$ is called an \emph{infra-nilmanifold}. If furthermore $E\subseteq N$, we say that the space $E\backslash N$ is a \emph{nilmanifold}.
\end{definition}

It is well known that infra-nilmanifolds are classified by their fundamental group that is almost-crystallographic~\cite{Au}. 
Every almost-crystallographic subgroup $E$ of the group $\aff{N}$ fits into an extension:
\begin{equation}\label{E:ses}
 1\to \Lambda \to E \to F \to 1,
\end{equation}
where $\Lambda=E\cap N$ is a uniform lattice in $N$, and $F$ is a finite subgroup of $C$ known as the \emph{holonomy group} of the corresponding infra-nilmanifold $E\backslash N$~\cite{Au}. 
Let $M$ be an infra-nilmanifold whose fundamental group $E$ is almost-crystallographic. 
Following~\cite[Page~788]{GPS},  we recall the construction of a faithful linear representation associated with the extension~\eqref{E:ses}. Suppose that the nilpotent lattice $\Lambda$ is of class $c+1$ \emph{i.e.}\ $\Gamma_c(\Lambda)\neq 1$ and $\Gamma_{c+1}(\Lambda)= 1$.
For $i=1,\ldots,c$, let $Z_i=\Gamma_i(\Lambda)\bigl/\Gamma_{i+1}(\Lambda)$ denote the factor groups of the lower central series $\{\Gamma_{i}(\Lambda)\}_{i=1}^{c+1}$ of $\Lambda$.
We assume that these quotients are torsion free, since this will be the case for the groups that we will study in the following. Thus, $Z_i\cong \Zz^{k_i}$ for all $1\leq i \leq c$, where $k_i>0$. The \emph{rank} or \emph{Hirsch number} or \emph{Hirsch length} of $\Lambda$ is equal to $\sum_{i=1}^{c}k_i$.
The action by conjugation of $E$ on $\Lambda$ induces an action of $E$ on $Z_{i}$, which factors through an action of the group $E/\Lambda$ (the holonomy group $F$), since $\Lambda$ acts trivially on $Z_{i}$.
This gives rise to a faithful representation $\theta_F\colon F\to \operatorname{\text{GL}}(n,\Zz)$ via the composition:
\begin{equation}\label{eq:theta}
 \theta_F\colon F\to \operatorname{\text{GL}}(k_1,\Zz)\times \cdots \times \operatorname{\text{GL}}(k_c,\Zz) \to \operatorname{\text{GL}}(n,\Zz),
\end{equation}
where $n$, the rank of $\Gamma$, is also equal to the dimension of~$N$.

In this work, we use part of the algebraic characterisation of almost-crystallographic groups given in~\cite[Theorem~3.1.3]{Dekimpe} as follows. 
Recall that a group $G$ is called \emph{polycyclic-by-finite} if it admits a subnormal series, that is, a sequence of subgroups $\{1\}= H_0 \leq H_1 \leq H_2\leq \cdots \leq H_n=G$ such that $H_i$ is normal in $H_{i+1}$, with $ \faktor{H_{i+1}}{H_i}$ a (infinite) cyclic or finite group for each $0\leq i\leq n-1$.

\begin{theorem}[{\cite[Theorem~3.1.3]{Dekimpe}}]\label{T:dekimpe}
Let $E$ be a polycyclic-by-finite group. Then $E$ is almost-crystallographic if and only if it has a  nilpotent subgroup, and possesses no non-trivial finite normal subgroups.
\end{theorem}

\section{Remarkable subgroups and quotients of \texorpdfstring{$\BB3$}{B3}} 
\label{S:presentation}

We start this section with a result on braid groups that holds for an arbitrary number of strands, but first we recall the following lemma that was proven in \cite{BDOS:2024}.

\begin{lemma}[{\cite[Lemma~3.11(2)]{BDOS:2024}}]
\label{L:8groups}

        Consider the following commutative diagrams of (vertical and horizontal) short exact sequences of groups in which every square is commutative
\begin{equation*}
\label{E:9groups2}
\xymatrix{
  &  & 1 \ar[d] & 1 \ar[d] & \\
1 \ar[r] & A \ar[r]^-{\iota_1} \ar[d]^-{\eta} & B \ar[r]^-{\rho_1} \ar[d]^-{\psi_1} & C \ar[r] \ar[d]^-{\psi_2} & 1\\
1 \ar[r] & D \ar[r]^-{\iota_2}  & E \ar[r]^-{\rho_2} \ar[d]^-{\phi_1} & F \ar[r] \ar[d]^-{\phi_2} & 1\\
 &  & G \ar[r]^-{\phi} \ar[d] & H  \ar[d] &  \\
  &   & 1  & 1  &
  }
\end{equation*}
Suppose that, for $i=1,2$, the homomorphisms $\iota_i$ and $\psi_i$ are inclusions. 
Then, $\eta$ is an isomorphism if and only if $\phi$ is.

\end{lemma}

In
this section we use several times the homomorphism $\rho_m$ from \eqref{eqn:rhom}.

\begin{proposition}
\label{P:isoquotient}
    Let $n\geq 3$ and $m\geq 3$  positive integers. 
  The groups $\faktor{B_n}{\rho_m^{-1}(Z(\rho_m(B_n)))}$ and $\faktor{\rho_m(B_n)}{Z(\rho_m(B_n))}$ are isomorphic.
\end{proposition}

\proof 
From the commutative diagram of vertical and horizontal short exact sequences
\begin{equation}
\label{E:diagram}
\xymatrix{
  &  & 1 \ar[d] & 1 \ar[d] & \\
1 \ar[r] & B_n[m] \ar[r] \ar@{=}[d] & \rho_m^{-1}(Z(\rho_m(B_n))) \ar[r]^-{\rho_m|} \ar[d] & Z(\rho_m(B_n)) \ar[r] \ar[d] & 1\\
1 \ar[r] & B_n[m] \ar[r]  & B_n \ar[r]^-{\rho_m} \ar[d] & \rho_m(B_n) \ar[r] \ar[d] & 1\\
 &  & \faktor{B_n}{\rho_m^{-1}(Z(\rho_m(B_n)))} \ar[r] \ar[d] & \faktor{\rho_m(B_n)}{Z(\rho_m(B_n))}  \ar[d] &  \\
  &   & 1  & 1  &
  }
\end{equation}
and Lemma~\ref{L:8groups} we get the result.

\endproof

\subsection{The level four braid group with three strands}\label{S:b34}

Let us recall that the congruence braid group $B_3[2]$ is isomorphic to the pure braid group on three strands~$P_3$~\cite{Arnold:1968}. 
We discuss a set of generators involved in the decomposition $P_3=\Zz\times F_2$. 
This is motivated by \cite[Proposition~8]{Goncalves-Guaschi:2004}, where the authors give a proof for a  decomposition of this kind, using the center $Z(P_n)=\Zz$, for any number $n$ of strands.

The pure braid group $P_3$ admits a presentation with generators $A_{1,2}, A_{1,3}, A_{2,3}$, see Subsection~\ref{Ss:pure} for the complete Artin presentation of pure braid groups. 

Consider the Fadell-Neuwirth short exact sequence
\[
\begin{CD}
1 @>>> F_2(A_{1,3}, A_{2,3}) @>>> P_3 @>{d_3}>> P_2 @>>> 1, 
\end{CD}
\]
where $d_3\colon P_3\to P_2$ is the homomorphism that ``geometrically forgets'' the last strand in $P_3$, $F_2(A_{1,3}, A_{2,3})$ is the rank 2 free group (freely generated by $A_{1,3}, A_{2,3}$) and $P_2$ is the pure braid group with 2 strands generated by~$A_{1,2}$. 
Consider the section of $d_3$ that sends $A_{1,2}\in P_2$ onto the full twist $\Delta_3^2\in P_3$. 
Since $Z(P_3)\cong \Zz$ is generated by the full twist we obtain the decomposition 
\begin{equation}
\label{eqn:p3}
P_3=\langle \Delta_3^2 \rangle \times F_2(A_{1,3}, A_{2,3}).
\end{equation}

Our first result in this subsection states that a similar result holds for the level four braid group with three strands.

\begin{proposition}\label{P:b34}
The group $B_3[4]$ admits a decomposition as $\Zz \times F_5$, where $F_5$ is a free group of rank 5 freely generated by $\{ A_{1,3}^2, A_{2,3}^2, A_{1,3} A_{2,3}^2 A_{1,3}^{-1}, [A_{2,3}, A_{1,3}], A_{1,3} [A_{2,3}, A_{1,3}] A_{1,3} \}$ and  $\Zz$ is generated by $\Delta_3^4$, the square of the full twist in $B_3$.
\end{proposition}

\proof
As recalled, the pure braid group $P_3$ admits a presentation with generators $A_{1,2}, A_{1,3},$  $A_{2,3}$. 
We shall also consider a presentation of $P_3$ coming from the decomposition given in \eqref{eqn:p3}. 
Let $\varphi\colon P_3\to \langle \Delta_3^2 \rangle \times F_2(A_{1,3}, A_{2,3})$ be the isomorphism defined by $\varphi(A_{1,2})=\Delta_3^2A_{2,3}^{-1}A_{1,3}^{-1}$, $\varphi(A_{1,3})=A_{1,3}$ and $\varphi(A_{2,3})=A_{2,3}$. 
Let $P_3^2$ denote the kernel of the composition $Ab_2\colon P_3\to \Zz\times \Zz \times \Zz \to \Zz/2\Zz \times \Zz/2\Zz \times \Zz/2\Zz$, obtained by abelianization mod~2. 
From \cite{Brendle-Margalit:2018} we know that $P_3^2$ is isomorphic to $B_3[4]$ the level four braid group with three strands. 

We obtain a presentation of $B_3[4]$ by computing a presentation of $\varphi(P_3^2)$. 
This is explained by the following commutative diagram of short exact sequences, where the vertical maps are all isomorphisms
\begin{equation*}
\xymatrix{
1 \ar[r] & P_3^2 \ar[r] \ar[d]^-{\varphi|} & P_3 \ar[r]^-{Ab_2} \ar[d]^-{\varphi} & \Zz/2\Zz \times \Zz/2\Zz \times \Zz/2\Zz \ar[r] \ar[d]^-{\overline{\varphi}} & 1\\
1 \ar[r] & Ker(\psi) \ar[r] & \langle \Delta_3^2 \rangle \times F_2(A_{1,3}, A_{2,3}) \ar[r]^-{\psi} & \Zz/2\Zz \times \Zz/2\Zz \times \Zz/2\Zz \ar[r] & 1}
\end{equation*}
and where $\psi$ is the abelianization mod~2 homomorphism for $\langle \Delta_3^2 \rangle \times F_2(A_{1,3}, A_{2,3})$, and the isomorphisms $\varphi|$ and $\overline{\varphi}$ are the restriction and projection of $\varphi$, respectively.

We note that $\psi=(\psi_1,\psi_2)$ where $\psi_1\colon \langle \Delta_3^2 \rangle\to \Zz/2\Zz$ and $\psi_2\colon F_2(A_{1,3}, A_{2,3})\to \Zz/2\Zz\times \Zz/2\Zz$ are the homomorphisms given by the respective abelianization mod~2. 
Hence, $Ker(\psi)=Ker(\psi_1)\times Ker(\psi_2)$. 
It is clear that $Ker(\psi_1)=\langle \Delta_3^4 \rangle$. It remains to characterize $Ker(\psi_2)$. Since $\Zz/2\Zz\times \Zz/2\Zz$ has 4 elements and $F_2(A_{1,3}, A_{2,3})$ a free group of rank 2, then $Ker(\psi_2)$ is free or rank 5 \cite[Theorem~2.10]{Magnus-Karrass-Solitar:book}. We can apply Reidemeister-Schreier method to compute the generators of $Ker(\psi_2)$. We choose the Schreier transversal to be $\{ 1, A_{1,3}, A_{2,3}, A_{1,3}A_{2,3} \}$. Then, by a direct calculation based on \cite[Theorem~2.9]{Magnus-Karrass-Solitar:book} we conclude that $Ker(\psi_2)$ is generated by $ A_{1,3}^2$, $A_{2,3}^2$, $A_{1,3} A_{2,3}^2 A_{1,3}^{-1}$, $[A_{2,3}, A_{1,3}]$, $A_{1,3} [A_{2,3}, A_{1,3}] A_{1,3} $.

\endproof

\begin{remark}
    A similar approach of Proposition~\ref{P:b34} for any number of strands will be one of the objects of a further work.
\end{remark}

From Proposition~\ref{P:b34} we obtain the following immediate consequence.

\begin{corollary}
    The consecutive quotients of the lower central series of the group $B_3[4]$ are free abelian groups such that
    \begin{enumerate}
        \item The abelianization of $B_3[4]$ is  a torsion free abelian group of rank $6$.
        
        \item For every $k\geq 2$, there is an isomorphism $\faktor{\Gamma_k(\BC3{4})}{\Gamma_{k+1}(\BC3{4})}\cong \faktor{\Gamma_k(F_5)}{\Gamma_{k+1}(F_5)}$, where $F_5$ is the free group of rank $5$.
    \end{enumerate}
\label{C:b34}
\end{corollary}

\proof
This follows from Proposition~\ref{P:b34} and basic properties of the lower central series of a group.
\endproof

We now move on to another useful consequence of Proposition~\ref{P:b34}. 

\begin{corollary}
\label{C:torsionfree}

\begin{enumerate}
    \item\label{i:first} The group $\rho_4^{-1}(Z(\rho_4(B_3)))$ is isomorphic to $\Zz \times F_5$, where $F_5$ is a free group of rank 5 freely generated by $\{ A_{1,3}^2, A_{2,3}^2, A_{1,3} A_{2,3}^2 A_{1,3}^{-1}, [A_{2,3}, A_{1,3}], A_{1,3} [A_{2,3}, A_{1,3}] A_{1,3} \}$ and  $\Zz$ is generated by the full twist $\Delta_3^2$ in $B_3$.

\item\label{i:second} For any $k\geq 2$ the equality $\Gamma_k(B_3[4]) = \Gamma_k\left(\rho_4^{-1}(Z(\rho_4(B_3)))\right)$ holds.

    \item\label{i:third} 
The group $\faktor{\rho_4^{-1}(Z(\rho_4(B_3)))}{[\BC3{4},\BC3{4}]}$ is isomorphic to $\Zz^6$, where one copy of $\Zz$ is generated by the class of the full twist $\Delta_3^2$.
\end{enumerate}
\end{corollary}

\proof 
Consider the short exact sequence 
\begin{equation}
\label{E:corb34}
\begin{CD}
1 @>>> \BC3{4} @>>> \rho_4^{-1}(Z(\rho_4(B_3))) @>{\rho_4|}>> Z(\rho_4(B_3)) @>>> 1. 
\end{CD}    
\end{equation}
\begin{enumerate}
    \item We recall that $Z(\rho_4(B_3))=\Zz/2\Zz$ is generated by the order two element $\rho_4(\Delta_3^2)$, see \cite[Lemma~2.3]{BDOS:2024}. 
We now consider the presentation of $\BC3{4} = \Zz[\Delta_3^4]\times F_5$ given in Proposition~\ref{P:b34}, where $\Zz[\Delta_3^4]$ is the cyclic group generated by $\Delta_3^4$. 
Then, from the method of presentation of extensions \cite[Chapter~10]{Johnson} we conclude that the group $\rho_4^{-1}(Z(\rho_4(B_3)))$  admits a presentation with generators given by $\Delta_3^2$ and the elements of the set $X=\{ A_{1,3}^2$, $A_{2,3}^2$, $A_{1,3} A_{2,3}^2 A_{1,3}^{-1}$, $[A_{2,3}, A_{1,3}]$, $A_{1,3} [A_{2,3}, A_{1,3}] A_{1,3}\}$ and relations $[\Delta_3^2,\, x]=1$ for all $x\in X$. From this computation, we get the point~\eqref{i:first} of this corollary.

\item From point~\eqref{i:first} of this corollary, we have that $\rho_4^{-1}(Z(\rho_4(B_3)))$ has a decomposition as $\Zz[\Delta_3^2] \times F_5[X]$ and from Proposition~\ref{P:b34} we know that $B_3[4]$ admits a decomposition as $\Zz[\Delta_3^4] \times F_5[X]$, where $F_5[X]$ is the free group freely generated by $X$. 
It is well known that the lower central series of a direct product of groups $G\times H$, $\Gamma_k(G\times H)$, is obtained by taking the direct product of the lower central series of $G$ and $H$, $\Gamma_k(G)\times \Gamma_k(H)$. 
Hence, we obtain point~\eqref{i:second} of this result.

\item This item is a consequence of the decomposition of $\rho_4^{-1}(Z(\rho_4(B_3)))$ given in \eqref{i:first} as the direct product $\Zz[\Delta_3^2] \times F_5[X]$, the equality $\Gamma_2(B_3[4]) = \Gamma_2\left(\rho_4^{-1}(Z(\rho_4(B_3)))\right)$ obtained in \eqref{i:second} and the elementary fact that the abelianization of the free group $F_5[X]$ is a free abelian group of rank 5.
\end{enumerate}
\endproof

\begin{remark}
\label{R:pres}

 Remark that the quotient group $\faktor{\rho_4(B_3)}{Z(\rho_4(B_3))}$ is isomorphic to $S_4$. 
 In fact, first we write a presentation of $\rho_4(B_3)$ using the decomposition (given in \cite{KordekMargalit:2022}) as a non-split extension of $S_3$ by $(\Zz/2\Zz)^3$ with generators: $A_{1,2}, A_{1,3}, A_{2,3}$ and $\sigma_1, \sigma_2$ and defining relations:
 \begin{enumerate}
 \item $\sigma_1^2=A_{1,2}$, $\sigma_2^2=A_{2,3}$, $\sigma_2\sigma_1^2\sigma_2^{-1}=A_{1,3}$, 
     \item $A_{i,j}^2=1$,
     \item $[A_{i,j}, A_{r,s}]=1$, for $1\leq i<j\leq 3$ and $1\leq r<s\leq 3$,
     \item $\sigma_1\sigma_2\sigma_1 = \sigma_2\sigma_1\sigma_2$,
     \item $\sigma_k A_{i,j} \sigma_k^{-1} = A_{\sigma_k(i), \sigma_k(j)}$, for $k=1,2$ and $1\leq i<j\leq 3$.
 \end{enumerate}
Then, we reduce this presentation to the following 
\begin{equation}
\label{E:pres}
    \langle \sigma_1,\, \sigma_2 \mid \sigma_1^4=1,\, \sigma_2^4=1,\, \sigma_1\sigma_2\sigma_1 = \sigma_2\sigma_1\sigma_2,\, [\sigma_1^2, \sigma_2^2]=1 \rangle.
\end{equation}
From this we conclude that the group $\rho_4(B_3)$ is isomorphic to $A_4\rtimes \Zz/4\Zz$,
its center is isomorphic to $\Zz/2\Zz$ and 
$\faktor{\rho_4(B_3)}{Z(\rho_4(B_3))}$ is isomorphic to $S_4$. 
This may be checked directly by hand or by using GAP-System \cite{GAP} with the presentation given in \eqref{E:pres}.
See Appendix~\ref{S:Appendix} for specific information about the algebraic structure of $\rho_4(B_3)$. 

\end{remark}

We finish this subsection by proving that the quotient group $\faktor{B_3}{[B_3[4], B_3[4]]}$ is a crystallographic group.

\begin{proposition}
    \label{P:b34cryst}
The group $\faktor{B_3}{[B_3[4], B_3[4]]}$ is crystallographic with dimension 6 and holonomy group~$S_4$.
\end{proposition}

\proof 
Consider the short exact sequence
    \[
\begin{CD}
1 @>>> \faktor{\BC3{4}}{[\BC3{4},\BC3{4}]} @>>> \faktor{B_3}{[\BC3{4},\BC3{4}]} @>{\overline{\rho_4}}>> \rho_4(B_3) @>>> 1.
\end{CD}
\]
Since $\Delta_3^2$ generates the center of $B_3$ and does not belong to $B_3[4]$, then from Theorem~\ref{T:33}  we know that the representation $\theta\colon \rho_4(B_3)\to \aut{\faktor{B_3[4]}{[B_3[4], B_3[4]]}}$, induced from the action by conjugation of $\faktor{B_3}{[\BC3{4},\BC3{4}]}$ on $\faktor{\BC3{4}}{[\BC3{4},\BC3{4}]}$, is not injective. 


Now, we compute the image of the element $A_{1,3}A_{2,3}^2A_{1,3}^{-1}\in B_3[4]$ (see Proposition~\ref{P:b34}) by the  automorphism $\theta(\sigma_2^2\sigma_1^2)$. Recall the action by conjugation of $B_3$ on $P_3$ described by using the Artin generators, see \cite[P.~5]{BDOS:2024}:
\[
\begin{aligned}
\sigma_1&:\ 
\begin{cases}
A_{1,2} \mapsto A_{1,2},\\[2pt]
A_{1,3} \mapsto A_{2,3},\\[2pt]
A_{2,3} \mapsto A_{2,3}^{-1}A_{1,3}A_{2,3}.
\end{cases}
\qquad
\sigma_2:\ 
\begin{cases}
A_{1,2} \mapsto A_{1,3},\\[2pt]
A_{1,3} \mapsto A_{1,3}^{-1}A_{1,2}A_{1,3},\\[2pt]
A_{2,3} \mapsto A_{2,3}.
\end{cases}
\end{aligned}
\]
Hence,
\begin{align*} 
\theta(\sigma_2^2\sigma_1^2)(A_{1,3}A_{2,3}^2A_{1,3}^{-1})  & = \theta(\sigma_2^2\sigma_1)\theta(\sigma_1)(A_{1,3}A_{2,3}^2A_{1,3}^{-1}) = \theta(\sigma_2^2\sigma_1)(A_{2,3}\cdot A_{2,3}^{-1}A_{1,3}^2A_{2,3}\cdot A_{1,3}^{-1})\\
& = \theta(\sigma_2^2)(A_{2,3}^2) = \theta(\sigma_2)(A_{2,3}^2) = A_{2,3}^2.
\end{align*} 
As a consequence, considering the classes of the six elements given in Proposition~\ref{P:b34} as a basis of the free abelian group  $\faktor{\BC3{4}}{[\BC3{4},\BC3{4}]}$, we conclude that $\sigma_2^2\sigma_1^2\notin Ker(\theta)$. 

We notice that $Ker(\theta)=Z(\rho_4(B_3))$. Indeed, from \cite[Lemma~2.3]{BDOS:2024} we know that $\rho_4(\Delta_3^2)$ is not trivial and has order 2 in $\rho_4(B_3)$. Then $Z(\rho_4(B_3))$ is contained in $Ker(\theta)$. Since $Ker(\theta)$ is a non-trivial normal subgroup of $\rho_4(B_3))$ then, from Theorem~\ref{T:app1}, it must be one of the groups $Z(\rho_4(B_3))$, $Z(\rho_4(B_3))\cdot \rho_4(B_3)^{(1)}$ or $Z(\rho_4(B_3))\cdot \rho_4(B_3)^{(2)}$. 
Since the element $\sigma_2^2\sigma_1^2\in \rho_4(B_3)^{(2)}$ does not belong to $Ker(\theta)$ then we conclude that $Ker(\theta)=Z(\rho_4(B_3))$.


From Remark~\ref{R:pres} follows that $\faktor{\rho_4(B_3)}{Z(\rho_4(B_3))} \cong S_4$. 
Consider now the following short exact sequence 
    \[
\begin{CD}
1 @>>> \faktor{\rho_4^{-1}(Z(\rho_4(B_3)))}{[\BC3{4},\BC3{4}]} @>>> \faktor{B_3}{[\BC3{4},\BC3{4}]} @>>> S_4  @>>> 1.
\end{CD}
\]
Since the group $\overline{\rho_4}^{-1}(Z(\rho_4(B_3)))=\faktor{\rho_4^{-1}(Z(\rho_4(B_3)))}{[\BC3{4},\BC3{4}]}$ is free abelian of dimension 6 (see Corollary~\ref{C:torsionfree}) then the result of this proposition follows from Theorem~\ref{T:34}.
\endproof

\subsection{On the group \texorpdfstring{$\BC3{p}$}{B3}}

We now explore some properties of the congruence subgroup $\BC3{p}$, where $p\geq 3$ is a prime. 
In the spirit of Proposition~\ref{P:b34}, we present here the following result. 

\begin{theorem}
    Let $p\geq 3$ be a prime and let $M = 1 + \frac{(p-1)(p)(p+1)}{12}$.
    The group $\BC3{p}$ is isomorphic to $\Zz\times F_M$, where $F_M$ is the free group of rank $M$. 
    \label{b3p_presentation}
\end{theorem}

\proof
Recall the Artin presentation for $\BB3$, with generators $\sig1, \sig2$ and  braid relation $\sig1 \sig2 \sig1 = \sig2 \sig1 \sig2$. 
Also, the special linear group $\mathrm{SL}_2(\mathbb{Z})$ has the following presentation (see \cite[Section 3.6.4]{Farb-Margalit:book})
\[ \langle a, b \mid a b a = b a b, (a b)^6=1 \rangle. \]
Consider the homomorphism $\BB3 \to \mathrm{SL}_2(\mathbb{Z})$ defined by $\sig1 \mapsto a$ and $\sig2 \mapsto b$. 
Let $\mathrm{SL}_2(\mathbb{Z})[p]$ be the kernel of $\mathrm{SL}_2(\mathbb{Z}) \to \mathrm{SL}_2(\mathbb{Z}/p\Zz)$. The image of the restriction of $\BB3 \to \mathrm{SL}_2(\mathbb{Z})$ to $\BC3{p}$ is $\mathrm{SL}_2(\mathbb{Z})[p]$ \cite[Theorem B]{bloomquist2023quotients}. 
Recall that $\Delta_3^2=(\sig1 \sig2)^3$ generates the center of $\BB3$, then the group $\Zz$ generated by $(\sig1 \sig2)^6$ is a normal subgroup of $B_3$ such that $\faktor{B_3}{\Zz}\cong SL_2(\Zz)$. 

Using the above information, we get the following diagram of (vertical and horizontal) short exact sequences of groups in which every square commutes
\begin{equation*}
\xymatrix{
  &  & 1 \ar[d] & 1 \ar[d] & \\
1 \ar[r] & Ker(\rho) \ar[r]^-{\iota} \ar[d]^-{\eta} & \BC3{p} \ar[r]^-{\rho} \ar[d]^-{\lambda} & \mathrm{SL}_2(\mathbb{Z})[p] \ar[r] \ar[d] & 1\\
1 \ar[r] & \Zz \ar[r]^-{\kappa}  & B_3 \ar[r] \ar[d] & SL_2(\Zz) \ar[r] \ar[d] & 1\\
 &  & \mathrm{SL}_2(\mathbb{Z}/p\Zz) \ar@{=}[r] \ar[d] & \mathrm{SL}_2(\mathbb{Z}/p\Zz) \ar[d] &  \\
  &   & 1  & 1  &
  }
\end{equation*}
From 
Lemma~\ref{L:8groups} 
we conclude that $\eta$ is an isomorphism. Even more, since $\iota$, $\kappa$ and $\lambda$ are inclusions, then $\eta$ is the identity isomorphism. Therefore $Ker(\rho)=\Zz$ is generated by $(\sig1 \sig2)^6$.

Since $p\geq 3$ is a prime, then $\mathrm{SL}_2(\mathbb{Z})[p]$ is a free group of rank $M = 1+ \frac{(p-1)(p)(p+1)}{12}$ \cite{Frasch:1933}. 
Let $t_0$ be the element $(\sig1 \sig2)^6$. Let  $t_1, t_2, \ldots, t_M$ denote elements of $\BC3{p}$ that map to generators of $\mathrm{SL}_2(\mathbb{Z})[p]$. Since $t_0$ commutes with every element of $\BC3{p}$ and the groups $\mathrm{SL}_2(\mathbb{Z})[p]$, $\mathbb{Z} = \langle t_0 \mid \, \rangle$ have no relations, by using the method of presentation of extensions \cite[Chapter~10]{Johnson} we have the following
\[ \BC3{p} \cong \langle t_0, t_1, t_2, \ldots, t_M \mid [t_0, t_i] = 1 \rangle. \]
From this presentation we conclude the result.
\endproof

\begin{remark}
Theorem~\ref{b3p_presentation} relies on the fact that for an odd prime $p$ the group $\mathrm{SL}_2(\Zz)[p]$ is free of rank $1 + \frac{(p-1)(p)(p+1)}{12}$ and $B_3[p]$ surjects onto $\mathrm{SL}_2(\Zz)[p]$. It is worth noting that for every positive integer $m$ such that $p|m$, we have that $\mathrm{SL}_2(\Zz)[m]$ is a free group of  finite rank as a finite index subgroup of $\mathrm{SL}_2(\Zz)[p]$. Moreover, since $\BC3{m}$ surjects onto $\mathrm{SL}_2(\Zz)[m]$ \cite[Theorem B]{bloomquist2023quotients}  the proof of Theorem~\ref{b3p_presentation} can be generalized to $\BC3{m}\cong \Zz \times F$, where $F$ is a free group of finite rank.
\end{remark}

As a corollary of Theorem~\ref{b3p_presentation} we get the following result.

\begin{corollary}
    Let $p\geq 3$ be a prime and let $M = 1 + \frac{(p-1)(p)(p+1)}{12}$. 
    The consecutive quotients of the lower central series of the group $\BC3{p}$ are free abelian groups such that
    \begin{enumerate}
        \item The abelianization of $\BC3{p}$ is  a torsion free abelian group of rank $M+1=2 + \frac{(p-1)(p)(p+1)}{12}$.
        \item for $k\geq 2$, $\faktor{\Gamma_k(\BC3{p})}{\Gamma_{k+1}(\BC3{p})}\cong \faktor{\Gamma_k(F_M)}{\Gamma_{k+1}(F_M)}$, where $F_M$ is the free group of rank $M$.
    \end{enumerate}
\label{b3p}
\end{corollary}

\proof
It follows from Theorem~\ref{b3p_presentation} and basic properties of the lower central series of a group. It is worth mentioning that the consecutive quotients of the  lower central series of a free group are free abelian groups, see \cite[ Corollary~5.12(iv)]{Magnus-Karrass-Solitar:book}.
\endproof

We also have the following useful result.

\begin{corollary}
\label{C:torsionfreeb3p}
Let $p\geq 3$ be a prime. 
\begin{enumerate}
    \item The group $\rho_p^{-1}(Z(\rho_p(B_3)))$ is isomorphic to $\Zz \times F_M$, where $F_M$ is a free group of rank $M = 1 + \frac{(p-1)(p)(p+1)}{12}$  and  $\Zz$ is generated by the full twist $\Delta_3^2$ in $B_3$.

\item For any $k\geq 2$ the equality $\Gamma_k(B_3[p]) = \Gamma_k\left(\rho_p^{-1}(Z(\rho_p(B_3)))\right)$ holds.

    \item  
The group $\faktor{\rho_p^{-1}(Z(\rho_p(B_3)))}{[\BC3{p},\BC3{p}]}$ is free abelian of rank $M+1$, where one copy of $\Zz$ is generated by the class of the full twist $\Delta_3^2$.
\end{enumerate}
\end{corollary}

\proof
This result is a consequence of Theorem~\ref{b3p_presentation} and its proof similar to the one of
Corollary~\ref{C:torsionfree}. 
\endproof

\subsection{The special case \texorpdfstring{$\BC3{3}$}{B3}} 

Recall that $\sigma_1, \sigma_2$ are the generators of $\BB3$. The group $\rho_3(B_3)$ is isomorphic to $\Sp{2}{\Zz/3\Zz}$ \cite[Theorem~1]{ACampo:1979}. Moreover, from \cite[Theorem~4.5]{Stylianakis:2018} $\BC3{3}$ is generated by the set
\begin{equation}\label{eq:gen_sty}
\{ \sigma_1^3,\, \sigma_2^3,\, \sigma_1 \sigma_2^3 \sigma_1^{-1},\, \sigma_1^{-1} \sigma_2^3 \sigma_1\}. 
\end{equation}
Consider the short exact sequence
\[
\begin{CD}
1 @>>> \faktor{\BC3{3}}{[\BC3{3},\BC3{3}]} @>>> \faktor{\BB3}{[\BC3{3},\BC3{3}]} @>>> \Sp{2}{\Zz/3\Zz} @>>> 1.
\end{CD}
\]

The group $\Sp{2}{\Zz/3\Zz}$ acts on $\faktor{\BC3{3}}{[\BC3{3},\BC3{3}]}$ as follows. Choose preimages of $\rho_3(\sigma_1), \rho_3(\sigma_2)$ in $\faktor{\BB3}{[\BC3{3},\BC3{3}]}$. Then, the elements given in \eqref{eq:gen_sty} act on elements of $\faktor{\BC3{3}}{[\BC3{3},\BC3{3}]}$ by conjugation. It is convenient to choose the cosets $\sigma_1 [\BC3{3},\BC3{3}], \sigma_2 [\BC3{3},\BC3{3}]$ as the preimages. In the proof of \cite[Theorem 4.5]{Stylianakis:2018} the action by conjugation of $\sig1, \sig2 \in \BB3$ on $\sigma^3_1, \sigma^3_2, \sigma_1 \sigma_2^3 \sigma_1^{-1}, \sigma_1^{-1} \sigma_2^3 \sigma_1 \in \BC3{3}$ is calculated. Therefore, the following actions are induced:


\[
\sigma_1\colon \left\{ \begin{matrix*}[l]
	\sigma^3_1 & \mapsto & \sigma_1^3 \bmod([B_3[3], B_3[3]]) \\
	\sigma^3_2  & \mapsto & \sigma_1 \sigma_2^3 \sigma_1^{-1} \bmod([B_3[3], B_3[3]]) \\
	\sigma_1 \sigma_2^3 \sigma_1^{-1}  & \mapsto & \sigma_1^{-1} \sigma_2^3 \sigma_1 \bmod([B_3[3], B_3[3]])\\
	\sigma^{-1}_1 \sigma_2^3 \sigma_1 & \mapsto & \sigma_2^3 \bmod([B_3[3], B_3[3]])
\end{matrix*}\right. 
\]

\[
\sigma_2\colon \left\{ \begin{matrix*}[l]
	\sigma^3_1 & \mapsto & \sigma_1^{-1} \sigma_2^3 \sigma_1 \bmod([B_3[3], B_3[3]]) \\
	\sigma^3_2  & \mapsto & \sigma^3_2 \bmod([B_3[3], B_3[3]]) \\
	\sigma_1 \sigma_2^3 \sigma_1^{-1} & \mapsto & \sigma^3_1 \bmod([B_3[3], B_3[3]])\\
	\sigma^{-1}_1 \sigma_2^3 \sigma_1 & \mapsto & \sigma_1 \sigma_2^3 \sigma_1^{-1} \bmod([B_3[3], B_3[3]])
\end{matrix*}\right. 
\]

Let $e_1=\sigma_1^3$, $e_2=\sigma_1\sigma_2^3\sigma_1^{-1}$, $e_3=\sigma_1^{-1}\sigma_2^3\sigma_1$, and $e_4=\sigma_2^3$ be elements in  $\faktor{\BC3{3}}{[\BC3{3}, \BC3{3}]}$. By Corollary \ref{b3p} the group $\faktor{\BC3{3}}{[\BC3{3},\BC3{3}]}$ is free abelian of rank 4. Since $\BC3{3}$ is generated by 4 elements \cite[Theorem 4.5]{Stylianakis:2018}, then the set $\{ e_1, e_2, e_3, e_4 \}$ is an ordered basis of $\faktor{\BC3{3}}{[\BC3{3}, \BC3{3}]}$. 


Therefore, the homomorphism $\Sp{2}{\Zz/3\Zz} \to \Aut\left(\faktor{\BC3{3}}{[\BC3{3},\BC3{3}]}\right)$ is determined by the action by conjugation of elements of $\faktor{\BB3}{[\BC3{3}, \BC3{3}]}$ on the ordered set $\{ e_1, e_2, e_3, e_4 \}$ of generators of the abelian group $\faktor{\BC3{3}}{[\BC3{3}, \BC3{3}]}$:
$$
\sigma_1: 
\begin{pmatrix}
1&0&0&0\\	
0&0&0&1\\
0&1&0&0\\	
0&0&1&0	
\end{pmatrix}
\quad
\textrm{ and }
\quad
\sigma_2: 
\begin{pmatrix}
0&1&0&0\\	
0&0&1&0\\
1&0&0&0\\	
0&0&0&1	
\end{pmatrix}.
$$

Recall from \cite{Goncalves-Guaschi-Ocampo:2017} that the group $\faktor{B_3}{[B_3[2], B_3[2]]}$ is not torsion free. 
In the next result we show that this does not hold when considering $B_3[3]$ instead of $B_3[2]$ in the quotient.

\begin{theorem}
The group $\faktor{\BB3}{[\BC3{3},\BC3{3}]}$ is torsion free.
\label{torsionfree_b3}
\end{theorem}
\proof 
Since, $\faktor{\BC3{3}}{[\BC3{3},\BC3{3}]} \cong \Zz^4$ and $\Sp2{\Zz/3\Zz}$ has 24 elements, from the short exact sequence
\begin{equation}\label{sesb33}
\begin{CD}
    1 @>>> \Zz^4 @>>> \faktor{\BB3}{[\BC3{3}, \BC3{3}]} @>\rho_3>> \Sp2{\Zz/3\Zz}  @>>> 1
\end{CD}
\end{equation}
we conclude that the only possible torsion in the middle group is 2, 3, 4 or 6. We will prove that $\faktor{\BB3}{[\BC3{3}, \BC3{3}]}$ does not have elements of order 2 or 3, and so it is torsion free.

We  note that $(\sigma_1^2\sigma_2)^2=\sigma_1^3\sigma_2^3$ in  $\faktor{\BB3}{[\BC3{3}, \BC3{3}]}$. The element $(\sigma_1^2\sigma_2)^2 = (\sig1 \sig2)^3$ is the full twist in $\BB3$ and generates its center. Also, $\rho_3((\sig1 \sig2)^3)$ is non-trivial and it has order 2, see \cite[Lemma 2.3]{BDOS:2024}. 

First, we verify that in $\faktor{\BB3}{[\BC3{3}, \BC3{3}]}$ the square of the element $(\sigma_1^2\sigma_2)^2$ is equal to $e_1e_3e_4e_2$:
\begin{align*}
  ((\sigma_1^2\sigma_2)^2)^2  & = ((\sigma_1\sigma_2)^3)^2\\
  & = \sigma_1\cdot \sigma_2\sigma_1\sigma_2\cdot \sigma_1\sigma_2\sigma_1\cdot \sigma_2\sigma_1\sigma_2\cdot \sigma_1\sigma_2\\
  & = \sigma_1\sigma_1\cdot  
  \sigma_2\sigma_1\sigma_2\cdot 
  \sigma_1\sigma_2\sigma_1\cdot 
  \sigma_2\sigma_1\sigma_1\sigma_2\\
  & = \sigma_1^3\sigma_2 \cdot \sigma_1\sigma_2\sigma_1 \cdot \sigma_2 
  \sigma_2\sigma_1\sigma_1\sigma_2\\
  & = \sigma_1^3\sigma_2^2\sigma_1\sigma_2^3\sigma_1^2\sigma_2\\
  & = \sigma_1^3\sigma_2^2\sigma_2\sigma_2^{-1}\sigma_1\sigma_2^3\sigma_1^{-1}\sigma_2\sigma_2^{-1}\sigma_1\sigma_1^2\sigma_2\\
& = \sigma_1^3\cdot \sigma_2^3 \cdot \sigma_2^{-1}\sigma_1\sigma_2^3\sigma_1^{-1}\sigma_2\cdot \sigma_2^{-1}\sigma_1^3\sigma_2\\
& = \sigma_1^3 \cdot \sigma_2^3 \sigma_2^{-1} \sigma_2^{-1} \sigma_1^{3} \sigma_2 \sigma_2 \cdot \sigma_1 \sigma_2^3 \sigma_1^{-1} \\
& = \sigma_1^3 \cdot \sigma_2^3 \sigma_2^{-2} \sigma_1^{3} \sigma_2^2 \cdot \sigma_1 \sigma_2^3 \sigma_1^{-1} \\
& = \sigma_1^3 \cdot \sigma_2^3 \sigma_2^{-3} \cdot \sigma_2 \sigma_1^{3} \sigma_2^{-1} \sigma_2^3 \cdot \sigma_1 \cdot \sigma_2^3 \sigma_1^{-1} \\
& = \sigma_1^3 \cdot \sigma_2 \sigma_1^{3} \sigma_1^{-1} \cdot \sigma_2^3 \cdot \sigma_1 \cdot \sigma_2^3 \sigma_1^{-1} \\
& = \sigma_1^3 \cdot \sigma_1^{-1} \sigma_2^{3} \sigma_2 \cdot \sigma_2^3 \cdot \sigma_1 \cdot \sigma_2^3 \sigma_1^{-1} \\
& = e_1 e_3 e_4 e_2 \in \faktor{\BC3{3}}{[\BC3{3}, \BC3{3}]}
\end{align*}
So, the element $(\sigma_1^2\sigma_2)^2$ projects onto the single element of order 2 in $\Sp2{\Zz/3\Zz}$.
The action of $(\sigma_1^2\sigma_2)^2$ on the ordered set $\{ e_1, e_2, e_3, e_4 \}$ of generators of the abelian group $\faktor{\BC3{3}}{[\BC3{3}, \BC3{3}]}$ is described by the identity matrix. 
Let $x_i\in \Zz$, for $i=1,2,3,4$. 
Then, if $(e_1^{x_1}e_2^{x_2}e_3^{x_3}e_4^{x_4} (\sigma_1^2\sigma_2)^2)^2$ is a trivial element in $\faktor{\BB3}{[\BC3{3}, \BC3{3}]}$ it implies that $e_1^{2x_1}e_2^{2x_2}e_3^{2x_3}e_4^{2x_4}e_1e_2e_3e_4$ is the identity element in $\faktor{\BC3{3}}{[\BC3{3}, \BC3{3}]}$. 
Since the equation $2x_j+1=0$ has no solutions in the set of integers  for $j=1,2,3,4$,  then $\faktor{\BB3}{[\BC3{3}, \BC3{3}]}$ does not have any element of order 2. 

Now, we prove that $\faktor{\BB3}{[\BC3{3}, \BC3{3}]}$ does not have elements of order 3. 
Remark that, since $\faktor{\BC3{3}}{[\BC3{3}, \BC3{3}]}$ is torsion free, an element of order 3 in $\faktor{\BB3}{[\BC3{3}, \BC3{3}]}$ projects onto an element of order 3 in $\Sp2{\Zz/3\Zz}$. 
Up to conjugacy, we may suppose that the order 3 elements in $\Sp2{\Zz/3\Zz}$ are either $\rho_3(\sigma_1)$ or $\rho_3(\sigma_1^2)$. 
Let $x_i\in \Zz$, for $i=1,2,3,4$. 
Then
\begin{align*}
    (e_1^{x_1}e_2^{x_2}e_3^{x_3}e_4^{x_4} \sigma_1)^3 & = e_1^{x_1}e_2^{x_2}e_3^{x_3}e_4^{x_4} \sigma_1e_1^{x_1}e_2^{x_2}e_3^{x_3}e_4^{x_4} \sigma_1e_1^{x_1}e_2^{x_2}e_3^{x_3}e_4^{x_4} \sigma_1\\
& = e_1^{x_1}e_2^{x_2}e_3^{x_3}e_4^{x_4} \sigma_1 e_1^{x_1}e_2^{x_2}e_3^{x_3}e_4^{x_4} \sigma_1^{-1} \sigma_1^2 e_1^{x_1}e_2^{x_2}e_3^{x_3}e_4^{x_4}\sigma_1^{-2} \sigma_1^3\\
& = e_1^{3x_1 + 1}e_2^{x_2+x_3+x_4}e_3^{x_2+x_3+x_4}e_4^{x_2+x_3+x_4}
\end{align*}
Since $3x_1 + 1=0$ has no solution in the set of integers, we conclude that $\faktor{\BB3}{[\BC3{3}, \BC3{3}]}$ does not have elements of order 3 that project onto the conjugacy class of $\rho_3(\sigma_1)$.
To prove that $\faktor{\BB3}{[\BC3{3}, \BC3{3}]}$ does not have elements of order 3 that project onto the conjugacy class of $\rho_3(\sigma_1^2)$ we use the same argument as above. The difference in this case is that $(\sigma_1^2)^3=e_1^2$ and so the problematic equation in the set of integers is, in this case, $3x_1+2=0$. 
\endproof

\section{Almost crystallographic groups}\label{S:almost}

We start this section noting that, from Theorem~\ref{T:dekimpe}, if $E$ is a crystallographic group then it does not have non-trivial finite normal subgroups. 
In \cite[Lemma~11(a)]{LIMAGONCALVES2019160} the authors proved that $\faktor{B_3[2]}{\Gamma_k(B_3[2])}$ is a torsion free  nilpotent group  of nilpotency class $k-1$ and its Hirsch length was computed. 
This motivates the following result, but now with congruence subgroups. That will be useful to show the connection between almost-crystallographic structures and congruence subgroups.

 \begin{lemma}\label{L:torsionfree}
 
Let $m$ be either 4 or an odd prime number and let $k\geq 2$. 
     \begin{enumerate}
\item\label{it:ptorsiona} The group $\faktor{B_3[m]}{\Gamma_k(B_3[m])}$ is a torsion free  nilpotent group  of nilpotency class $k-1$ with Hirsch length 
equal to 
\begin{equation*}
\sum_{q=1}^{k-1}\left( \frac{1}{q}\sum_{d\mid q}\mu(d)M^{q/d} \right) + 1, 
\end{equation*}
where $\mu$  is the M\"obius function and $M$ is equal to $5$ if $m=4$ or $M = 1 + \frac{(m-1)(m)(m+1)}{12}$ if $m$ is an odd prime.

\item\label{it:ptorsionb} The quotient group $\faktor{B_3}{\Gamma_k(B_3[m])}$ does not have non-trivial normal finite subgroups. 
\end{enumerate}
 \end{lemma}

\proof 
Let $m$ be either 4 or an odd prime number and let $k\geq 2$.  
     \begin{enumerate}
         \item  
         Let $M$ be equal to $5$ if $m=4$ or $M = 1 + \frac{(m-1)(m)(m+1)}{12}$ if $m$ is an odd prime. Let $F_M$ denote the free group of rank $M$. 
From the Witt's Formulae, see \cite[Theorem~5.11]{Magnus-Karrass-Solitar:book}, the rank of the free abelian group $\faktor{\Gamma_k(F_M)}{\Gamma_{k+1}(F_M)}$ is equal to $\frac{1}{k}\sum_{d\mid k}\mu(d)M^{k/d}$, where $\mu$  is the M\"obius function. 
From this and Corollary~\ref{C:b34} if $m=4$, or Corollary~\ref{b3p} if $m$ is an odd prime, it follows that the  Hirsch length of the  nilpotent group $\faktor{B_3[m]}{\Gamma_k(B_3[m])}$ of nilpotency class $k-1$ is equal to 
$\sum_{q=1}^{k-1}\left( \frac{1}{q}\sum_{d\mid q}\mu(d)M^{q/d} \right) + 1.$    
The fact that $\faktor{B_3[m]}{\Gamma_k(B_3[m])}$ is torsion free follows by induction using the short exact sequence
\begin{equation*}
1 \to \faktor{\Gamma_{j-1}(B_3[m])}{\Gamma_j(B_3[m])} \to \faktor{B_3[m]}{\Gamma_j(B_3[m])} \to \faktor{B_3[m]}{\Gamma_{j-1}(B_3[m])} \to 1
\end{equation*}
since the abelianisation group $\faktor{B_3[m]}{\Gamma_{2}(B_3[m])}$ and the consecutive quotients of the lower central series $\faktor{\Gamma_{j-1}(B_3[m])}{\Gamma_j(B_3[m])}$ are torsion free groups.

         \item Let $j\geq 2$. 
         The quotient group $\faktor{\Gamma_{j-1}(B_3[m])}{\Gamma_j(B_3[m])}$ is torsion-free for any $j$, see
         Corollary~\ref{C:b34} if $m=4$ or Corollary~\ref{b3p} if $m$ is an odd prime.
The group $\faktor{B_3}{\Gamma_{2}(B_3[m])}$ does not have non-trivial normal finite subgroups since it is a crystallographic group.  This follows from Proposition~\ref{P:b34cryst} if $m=4$, 
or from \cite[Theorem~3.6]{BDOS:2024} if $m$ is an odd prime.
         We get the proof of this item  by induction on $j$, since for all $j\geq 2$, we have a central extension of the form:
\begin{equation}\label{E:gammaquot}
1 \to \faktor{\Gamma_{j-1}(B_3[m])}{\Gamma_j(B_3[m])} \to \faktor{B_3}{\Gamma_j(B_3[m])} \to \faktor{B_3}{\Gamma_{j-1}(B_3[m])} \to 1.
\end{equation}
     \end{enumerate}
\endproof

We are now ready to state the main result of this section.

\begin{theorem}
Let $m$ be either 4 or an odd prime and let $k\geq 2$ be a positive integer. Let also $M=5$ if $m=4$ or $M = 1 + \frac{(m-1)(m)(m+1)}{12}$ if $m$ is an odd prime. Finally, denote by $G_m$ the group $S_4$ if $m=4$ or the quotient $\faktor{\Sp{2}{\Zz/m\Zz}}{Z(\Sp{2}{\Zz/m\Zz})}$ if $m$ is an odd prime. There is a short exact sequence
\[
\begin{CD}
1 \to  \faktor{\rho_m^{-1}(Z(\rho_m(B_3)))}{\Gamma_k(\BC3{m})} \to  \faktor{B_3}{\Gamma_k(\BC3{m})} \to  \faktor{\rho_m(B_3)}{Z(\rho_m(B_3))} \to  1,
\end{CD}
\] 
where the middle
group $\faktor{B_3}{\Gamma_k(\BC3{m})}$ is an almost-crystallographic group whose holonomy group is $G_m$ and whose dimension is equal to
\[
\sum_{q=1}^{k-1}\left( \frac{1}{q}\sum_{d\mid q}\mu(d)M^{q/d} \right) + 1,
\]
where $\mu$  is the M\"obius function. 
\label{T:alm_cryst}
\end{theorem}

\proof

Let $m$ be either 4 or an odd prime and let $k\geq 2$ be a positive integer. 
The short exact sequence
$
\begin{CD}
1 \to  \rho_m^{-1}(Z(\rho_m(B_3))) \to  B_3 \to  \faktor{B_3}{\rho_m^{-1}(Z(\rho_m(B_3)))} \to  1
\end{CD}
$ 
induces the following short exact sequence
\[
\begin{CD}
1 \to \faktor{\rho_m^{-1}(Z(\rho_m(B_3)))}{\Gamma_k(\BC3{m})} \to  \faktor{B_3}{\Gamma_k(\BC3{m})} \to  \faktor{B_3}{\rho_m^{-1}(Z(\rho_m(B_3)))} \to  1
\end{CD}
\]
From Proposition~\ref{P:isoquotient} we have the isomorphism 
$\faktor{\rho_m(B_3)}{Z(\rho_m(B_3))} \cong \faktor{B_3}{\rho_m^{-1}(Z(\rho_m(B_3)))}$ which leads us to obtain the short exact sequence of the statement. 
 Recall that $\faktor{\rho_m(B_3)}{Z(\rho_m(B_3))}$ is isomorphic to $S_4$ if $m=4$ (see Remark~\ref{R:pres}) or to $\faktor{\Sp{2}{\Zz/m\Zz}}{Z(\Sp{2}{\Zz/m\Zz})}$ if $m$ is an odd prime.

From the above, Theorem~\ref{T:dekimpe} and the second item of  Lemma~\ref{L:torsionfree} we prove that the group $\faktor{B_3}{\Gamma_k(\BC3{m})}$ is an almost-crystallographic group that fits in the extension
\[
\begin{CD}
1 \to  \faktor{\rho_m^{-1}(Z(\rho_m(B_3)))}{\Gamma_k(\BC3{m})} \to  \faktor{B_3}{\Gamma_k(\BC3{m})} \to \faktor{\rho_m(B_3)}{Z(\rho_m(B_3))} \to  1.
\end{CD}
\]

We note that, from the first item of Lemma~\ref{L:torsionfree} and from Corollaries~\ref{C:torsionfree}
and~\ref{C:torsionfreeb3p}, the group $\faktor{\rho_m^{-1}(Z(\rho_m(B_3)))}{\Gamma_k(\BC3{m})}$ is a torsion free  nilpotent group  of nilpotency class $k-1$ with Hirsch length equal to 
\begin{equation*}
\sum_{q=1}^{k-1}\left( \frac{1}{q}\sum_{d\mid q}\mu(d)M^{q/d} \right) + 1, 
\end{equation*}
where $\mu$  is the M\"obius function and $M$ is equal to $5$ if $m=4$ or $M = 1 + \frac{(m-1)(m)(m+1)}{12}$ if $m$ is an odd prime.
\endproof

Notice that $A_4$ is isomorphic to $\faktor{\Sp2{\Zz/3\Zz}}{Z(\Sp2{\Zz/3\Zz})}$ \cite[Exercise 5.8]{Fulton-Harris:1991}. We have the following result for the special case $m=3$.

\begin{corollary}
Let $k\geq 2$. 
Then $\faktor{\BB3}{\Gamma_k(\BC3{3})}$ is almost-Bieberbach with holonomy group $A_4$. 
In particular, the group $\faktor{\BB3}{[\BC3{3},\BC3{3}]}$ is Bieberbach of dimension 4 and holonomy group $A_4$.
\end{corollary}

\proof
Let $k\geq 2$. 
From Theorem~\ref{T:alm_cryst} the group $\faktor{\BB3}{\Gamma_k(\BC3{3})}$ is almost-crystallographic with holonomy group $A_4$.

From \cite[Theorem~3.6]{BDOS:2024} the group $\faktor{\BB3}{[\BC3{3},\BC3{3}]}$ is crystallographic of dimension 4 and holonomy group $A_4$ and from Theorem~\ref{torsionfree_b3} it is torsion free. Hence it is a Bieberbach group.

Now we use again the fact that $\faktor{\BB3}{[\BC3{3},\BC3{3}]}$ is torsion free for the base step ($j=3$) in an induction process. 
The general case of the induction follows from the extension given in \eqref{E:gammaquot} since the kernel of that extension is torsion free, see Corollary~\ref{b3p}.
\endproof

\begin{remark}
    
For $m=3$ we note the following: $\Sp2{\Zz/3\Zz}$ has 7 conjugacy classes of elements. The only non-trivial class which is mapped to the identity in $\mathrm{Aut}(\mathbb{\Zz}^4)$ is the class of $\rho_3((\sigma_1^2\sigma_2)^2)$. 
This element has order 2 and generates the center of $\Sp2{\Zz/3\Zz}$. Furthermore, it belongs to the subgroups of $\Sp2{\Zz/3\Zz}$ with the exception of $\Zz/3\Zz$.
Hence, we may conclude that $\faktor{\rho_3^{-1}(\Zz/3\Zz)}{[\BC3{3}, \BC3{3}]}$ is a crystallographic group (in fact Bieberbach) of dimension 4 and holonomy group $\Zz/3\Zz$.

Since $\faktor{\rho_3^{-1}(\Zz/3\Zz)}{[\BC3{3}, \BC3{3}]}$ is torsion free, it is a Bieberbach group. Then, from Bieberbach theorems (see \cite[Section~2.1]{Dekimpe}), there is a 4-dimensional compact Riemannian flat manifold $\chi$ with holonomy group $\Zz/3\Zz$ that has fundamental group $\faktor{\rho_3^{-1}(\Zz/3\Zz)}{[\BC3{3}, \BC3{3}]}$.
Motivated by results about the holonomy representation of Bieberbach subgroups of the Artin braid group quotient $\faktor{B_n}{[P_n, P_n]}$ whose holonomy group is a 2-group obtained in \cite[Section~3]{OR}, applying a result of Porteous \cite[Theorem~7.1]{Porteous} is it not difficult to prove that $\chi$ admits Anosov diffeomorphism.
\end{remark}


\appendix

\section{On the finite group $\rho_4(B_3)$}\label{S:Appendix}


From the definitions of congruence subgroups given in Subsection~\ref{Ss:cong}, taking the specific cases $n=3$ and $m=4$ in  equation~\eqref{eqn:rhom}, we have the symplectic representation $\rho_4\colon \BB3 \rightarrow \Sp{2}{\Zz/4\Zz}$. In this situation, the kernel of $\rho_4$ is $\BC3{4}$, the level $4$ congruence subgroup of $\BB3$. 
In this appendix, we study the algebraic structure of the finite group $\rho_4(B_3)$, mentioned in Remark~\ref{R:pres}. The properties given here are useful in Subsection~\ref{S:b34}.

We start recalling the presentation of $\rho_4(B_3)$ using the decomposition (given in \cite{KordekMargalit:2022}) as a non-split extension of $S_3$ by $(\Zz/2\Zz)^3$ with generators: $A_{1,2}, A_{1,3}, A_{2,3}$ and $\sigma_1, \sigma_2$ and defining relations:
 \begin{enumerate}
 \item $\sigma_1^2=A_{1,2}$, $\sigma_2^2=A_{2,3}$, $\sigma_2\sigma_1^2\sigma_2^{-1}=A_{1,3}$, 
     \item $A_{i,j}^2=1$,
     \item $[A_{i,j}, A_{r,s}]=1$, for $1\leq i<j\leq 3$ and $1\leq r<s\leq 3$,
     \item $\sigma_1\sigma_2\sigma_1 = \sigma_2\sigma_1\sigma_2$,
     \item $\sigma_k A_{i,j} \sigma_k^{-1} = A_{\sigma_k(i), \sigma_k(j)}$, for $k=1,2$ and $1\leq i<j\leq 3$.
 \end{enumerate}
We simplify it to the following 
\begin{equation}
\label{E:Apres}
    \rho_4(B_3) = \langle \sigma_1,\, \sigma_2 \mid \sigma_1^4=1,\, \sigma_2^4=1,\, \sigma_1\sigma_2\sigma_1 = \sigma_2\sigma_1\sigma_2,\, [\sigma_1^2, \sigma_2^2]=1 \rangle.
\end{equation}

\begin{remark}
 Notice that $\rho_4(B_3)$ is the small group of order 48 and ID 30 among the groups of order 48 in \emph{GAP's SmallGroup library} \cite{GAP}.
 
\begin{enumerate}
    \item It is easy to verify it in GAP by using the command \emph{IdGroup(G)}, obtaining as answer $[48,\, 30]$, defining the group $G$ in GAP using a presentation of it, e.g. \eqref{E:Apres}. 

Furthermore, the ID of the groups $\Sp{2}{\Zz/4\Zz}$ and $SL_2(\Zz/4\Zz)$ in GAP is exactly $[48,\, 30]$. Hence, we conclude that $\rho_4(B_3)$ is isomorphic to both $\Sp{2}{\Zz/4\Zz}$ and $SL_2(\Zz/4\Zz)$. 

\item The group $SL_2(\Zz/4\Zz)$ has been studied, see for instance \url{https://people.maths.bris.ac.uk/~matyd/GroupNames/1/A4sC4.html} where several properties of this group are exhibited. 
Although $SL_2(\Zz/4\Zz)$ is well understood, we would like to analyse $\rho_4(B_3)$ by using the generators $\sigma_1$, $\sigma_2$ of \eqref{E:Apres}, since it allows us to make the connection with braid groups that is useful in this paper.

\end{enumerate}
\end{remark}

In this appendix, we use the well-known semidirect product presentation
\begin{equation}
\label{E:Apres2}
A_4 \;=\;
\big\langle \alpha,\, r,\, s \;\big|\;
\alpha^3=1,\; r^2=1,\; s^2=1,\; rs=sr,\;
\alpha r\alpha^{-1}=s,\;
\alpha s\alpha^{-1}=rs
\big\rangle,
\end{equation}
corresponding to the decomposition \(A_4 \cong V \rtimes \Zz/3\Zz\),
where \(\langle r,s\rangle \cong V\) is the Klein four group.

Now, we modify \eqref{E:Apres} to show that the group $\rho_4(B_3)$ is isomorphic to $A_4\rtimes \Zz/4\Zz$, by giving a presentation that reveals this fact. First, we show that in the group $\rho_4(B_3)$ there are other elements with nice relations. 

\begin{proposition}\label{P:app1}
Let $\alpha = \sigma_2\sigma_1^{-1}$, $r=\sigma_2^2\sigma_1^2$ and $s=\sigma_1(\sigma_2^2\sigma_1^2)\sigma_1^{-1}$ in $\rho_4(B_3)$. Then, the following relations hold in $\rho_4(B_3)$: 
\begin{multicols}{3}
\begin{enumerate}
    \item $\alpha^3=1$
    \item $\alpha r\alpha^{-1} = s$
    \item $\alpha s\alpha^{-1} = rs$
    \item $r^2=1$
    \item $s^2=1$
    \item $rs=sr$
    \item $\sigma_1\alpha\sigma_1^{-1} = \alpha^{-1}r$
    \item $\sigma_1 r\sigma_1^{-1} = s$
    \item $\sigma_1 s\sigma_1^{-1} = r$
\end{enumerate}
\end{multicols}
\end{proposition}

\proof 
Let $\alpha = \sigma_2\sigma_1^{-1}$, $r=\sigma_2^2\sigma_1^2$ and $s=\sigma_1(\sigma_2^2\sigma_1^2)\sigma_1^{-1}$ in $\rho_4(B_3)$. 
We prove the following equalities by using the relations of the group $\rho_4(B_3)$ given in \eqref{E:Apres}.

\begin{enumerate}
    \item   \begin{align*}
            \alpha^3  & = \sigma_2\sigma_1^{-1}\sigma_2\sigma_1^{-1}\sigma_2\sigma_1^{-1} = \sigma_2\sigma_1^{3}\sigma_2\sigma_1^{-1}\sigma_2\sigma_1^{-1}  = \sigma_2\sigma_1^{2}\cdot \sigma_2^{-1}\sigma_1\sigma_2\cdot \sigma_2\sigma_1^{-1} \\
                    & = \sigma_2\sigma_1^{2}\sigma_2^{-1}\cdot\sigma_2^{-1} \sigma_1^2\sigma_2 = \sigma_2\sigma_1^{2}\sigma_2^{-2} \sigma_1^2\sigma_2 = \sigma_2\sigma_2^{-2}\sigma_2 = 1.
    \end{align*}

    \item   \begin{align*} 
            \alpha r\alpha^{-1} & = \sigma_2\sigma_1^{-1} \cdot \sigma_2^2\sigma_1^2\cdot \sigma_1\sigma_2^{-1} = \sigma_2\sigma_1 \sigma_2^2 \sigma_1\sigma_2^{-1} = \sigma_1\sigma_2 \sigma_1\cdot \sigma_1^{-1}\sigma_2\sigma_1   \\
            & = \sigma_1\sigma_2^2 \sigma_1 = \sigma_1(\sigma_2^2\sigma_1^2) \sigma_1^{-1} = s. 
            \end{align*}
    
    \item   \begin{align*} 
            \alpha s\alpha^{-1} & = \sigma_2\sigma_1^{-1} \cdot \sigma_1(\sigma_2^2\sigma_1^2)\sigma_1^{-1}\cdot \sigma_1\sigma_2^{-1} =  \sigma_2^2 \sigma_2\sigma_1^2\sigma_2^{-1} =  \sigma_2^2 \cdot \sigma_1^{-1}\sigma_2^2\sigma_1\\
            & =  \sigma_2^2 \sigma_1^3\sigma_2^2\sigma_1 =  \sigma_2^2 \sigma_1^2\cdot \sigma_1(\sigma_2^2\sigma_1^2)\sigma_1^{-1} = rs.
            \end{align*}

    \item \begin{align*} 
            r^2 & =  \sigma_2^2\sigma_1^2\cdot \sigma_2^2\sigma_1^2 = \sigma_2^2\sigma_2^2\cdot \sigma_1^2\sigma_1^2 = 1.
            \end{align*}
   
    \item \begin{align*} 
            s^2 & = \sigma_1(\sigma_2^2\sigma_1^2)\sigma_1^{-1}\cdot \sigma_1(\sigma_2^2\sigma_1^2)\sigma_1^{-1} = \sigma_1(\sigma_2^2\sigma_1^2\sigma_2^2\sigma_1^2)\sigma_1^{-1} = 1.
            \end{align*}

    \item Since $r^2=s^2=1$, then to prove $rs=sr$ we prove $(rs)^2=1$.
    \begin{align*} 
            (rs)^2 & = \sigma_2^2\sigma_1^2\cdot \sigma_1(\sigma_2^2\sigma_1^2)\sigma_1^{-1}\cdot  \sigma_2^2\sigma_1^2\cdot \sigma_1(\sigma_2^2\sigma_1^2)\sigma_1^{-1} = \sigma_2^2\sigma_1^3\sigma_2^2\sigma_1\sigma_2^2\sigma_1^3\sigma_2^2\sigma_1 \\
            & = \sigma_2^2\cdot \sigma_2\sigma_1^2\sigma_2^{-1}\cdot \sigma_2^2\cdot \sigma_2\sigma_1^2\sigma_2^{-1} = \sigma_2^2\sigma_2\sigma_1^2\sigma_2^2\sigma_1^2\sigma_2^{-1} = 1.
            \end{align*}
            
        \item \begin{align*} 
            \sigma_1\alpha\sigma_1^{-1} & = \sigma_1 (\sigma_2\sigma_1^{-1}) \sigma_1^{-1} = \sigma_1 \cdot \sigma_2^{-1}\sigma_2\cdot \sigma_2\sigma_1^2 = (\sigma_2\sigma_1^{-1})^{-1}\cdot \sigma_2^2\sigma_1^2 = \alpha^{-1}r.
            \end{align*}

        \item \begin{align*} 
            \sigma_1 r \sigma_1^{-1} & = \sigma_1(\sigma_2^2\sigma_1^2)\sigma_1^{-1} = s.
            \end{align*}
        \item \begin{align*} 
            \sigma_1 s \sigma_1^{-1} & = \sigma_1(\sigma_1\sigma_2^2\sigma_1^2\sigma_1^{-1})\sigma_1^{-1} = \sigma_2^2\sigma_1^2 = r.
            \end{align*}
\end{enumerate}

\endproof

\begin{remark}\label{R:app1}
Consider the elements $\alpha$, $r$ and $s$ of Proposition~\ref{P:app1}. 
\begin{enumerate}
\item Notice that in $\rho_4(B_3)$ the equality $\sigma_2^4=\sigma_2^2\sigma_1^2\cdot \sigma_2^2\sigma_1^2$ holds, and so the relation $\sigma_2^4=1$ is equivalent to $r^2=1$. Also, the braid relation $\sigma_1\sigma_2\sigma_1 = \sigma_2\sigma_1\sigma_2$ becomes $\sigma_1\alpha\sigma_1^{-1} = \alpha^{-1}r$. 

\item The image of the full-twist $\rho_4(\Delta_3^2)=\rho_4((\sigma_1\sigma_2)^3)$ is equal to $\sigma_1\alpha\sigma_1^{-1}\alpha\sigma_1^2$. 
Indeed:
\begin{align*} 
            \rho_4(\Delta_3^2) & = \sigma_1\sigma_2 \sigma_1\sigma_2 \sigma_1\sigma_2 = \sigma_1\sigma_2 \sigma_1^2 \sigma_2\sigma_1 = \sigma_1\sigma_2^{-1} \sigma_2^2 \sigma_1^2 \sigma_2\sigma_1^{-1}\sigma_1^{-2} \\
            & = \alpha^{-1} r \alpha \sigma_1^2 = \sigma_1\alpha\sigma_1^{-1}\alpha\sigma_1^2.
            \end{align*}
    In the last equality we used Proposition~\ref{P:app1}(7).

\item The subgroup of $\rho_4(B_3)$ generated by $\alpha$, $r$ and $s$ with the defining relations (1) to (6) is isomorphic to the alternating group $A_4=V\rtimes \Zz/3\Zz$, where $V$ is the Klein four group generated by $\{r,\, s\}$ and $\Zz/3\Zz$ is generated by $\alpha$, see \eqref{E:Apres2}. 

\end{enumerate}
\end{remark}

Recall that a group $G$ is \emph{solvable} is its derived series eventually reaches the trivial subgroup of $G$. The least $n$ such that $G^{(n)} = 1$ is called the \emph{derived length} of the solvable group $G$. 
Also recall that  subgroup $H$ of a group $G$ is called a \emph{characteristic subgroup} if for every automorphism $\phi$ of $G$, one has $\phi(H)=H$.

\begin{theorem}\label{T:app1}
Let $\alpha = \sigma_2\sigma_1^{-1}$, $r=\sigma_2^2\sigma_1^2$ and $s=\sigma_1(\sigma_2^2\sigma_1^2)\sigma_1^{-1}$ in $\rho_4(B_3)$. 
The group $\rho_4(B_3)$ admits the presentation
\begin{equation}
\label{E:Apres3}
\bigg\langle \sigma_1,\, \alpha,\, r,\, s \ \bigg\vert \ 
\begin{matrix}
\alpha^3=1,\; r^2=1,\; s^2=1,\; rs=sr,\; \alpha r\alpha^{-1}=s,\; 
\alpha s\alpha^{-1}=rs,\; \\ \sigma_1^4=1,\; 
\sigma_1\alpha\sigma_1^{-1} = \alpha^{-1}r,\;
\sigma_1 r\sigma_1^{-1} = s\; 
\sigma_1 s\sigma_1^{-1} = r
\end{matrix}
\ \bigg\rangle.
\end{equation}
Hence, it is isomorphic to $A_4\rtimes \Zz/4\Zz$, where $A_4=V\rtimes \Zz/3\Zz$ is the group generated by $\{\alpha,\, r,\, s\}$ and $\Zz/4\Zz$ is generated by $\sigma_1$.  
Furthermore,
\begin{enumerate}
    \item The group $\rho_4(B_3)$ is solvable, of derived length 3, with derived subgroup $\rho_4(B_3)^{(1)}=A_4$ and second derived subgroup the Klein four group $\rho_4(B_3)^{(2)}=V$.
    \item The center $Z(\rho_4(B_3))$ has order 2, it is generated by $\sigma_1\alpha\sigma_1^{-1}\alpha\sigma_1^2$.
    \item The group $\rho_4(B_3)$ has 5 non-trivial characteristic subgroups, which are pairwise different. They are 
    \begin{multicols}{2}
\begin{enumerate}
    \item $Z(\rho_4(B_3))=\Zz/2\Zz$, 
    \item $\rho_4(B_3)^{(1)}=A_4$, 
    \item $\rho_4(B_3)^{(2)}=V$, 
    \item $Z(\rho_4(B_3))\cdot \rho_4(B_3)^{(1)} \cong \Zz/2\Zz \times A_4$, and 
    \item $Z(\rho_4(B_3))\cdot \rho_4(B_3)^{(2)}\cong \Zz/2\Zz \times V$. 
\end{enumerate}
\end{multicols}
\end{enumerate}

\end{theorem}

\proof
Let $\alpha = \sigma_2\sigma_1^{-1}$, $r=\sigma_2^2\sigma_1^2$ and $s=\sigma_1(\sigma_2^2\sigma_1^2)\sigma_1^{-1}$ in $\rho_4(B_3)$. 
From Proposition~\ref{P:app1}, Remark~\ref{R:app1} and the presentation of the group $\rho_4(B_3)$ given in \eqref{E:Apres} we obtain the presentation \eqref{E:Apres3} of $\rho_4(B_3)$, where the action of $\Zz/4\Zz$ on $A_4$ is described in Proposition~\ref{P:app1} items (7) to (9).

From the presentation \eqref{E:Apres3} is clear that $\rho_4(B_3)$ is isomorphic to $A_4\rtimes \Zz/4\Zz$, where $A_4=V\rtimes \Zz/3\Zz$ is the group generated by $\{\alpha,\, r,\, s\}$ and $\Zz/4\Zz$ is generated by $\sigma_1$.

\begin{enumerate}
    \item Follows from \eqref{E:Apres3} that the abelianization of the group $\rho_4(B_3)$ is isomorphic to $\Zz/4\Zz$, generated by the class of $\sigma_1$ and that the commutator subgroup is $A_4$. It is well known that the commutator subgroup of $A_4$ is the Klein four group $V$.  Therefore, the group $\rho_4(B_3)$ is solvable, of derived length 3, with derived subgroup $\rho_4(B_3)^{(1)}=A_4$ and second derived subgroup the Klein four group $\rho_4(B_3)^{(2)}=V$.

    \item Recall from \cite[Lemma~2.3]{BDOS:2024} that $\rho_4(\Delta_3^2)$ is not trivial and has order 2 in $\rho_4(B_3)$. Since $Z(B_n)=\langle \Delta_3^2 \rangle$ then $Z(\rho_4(B_3)) = \langle \rho_4(\Delta_3^2) \rangle$. From Remark~\ref{R:app1}, the element $\rho_4(\Delta_3^2)$ has a description using the generators given in \eqref{E:Apres3} as $\sigma_1\alpha\sigma_1^{-1}\alpha\sigma_1^2$.

    \item To verify this item, we use the following code in GAP-System \cite{GAP}, to list all normal subgroups of $\rho_4(B_3)$:
    \begin{lstlisting}[language=GAP]
f2:=FreeGroup("u","v");;
AssignGeneratorVariables(f2);;
rel:=ParseRelators([u,v], "u^4=1, v^4=1, uvu=vuv, u^2v^2=v^2u^2");;
rho4B3:=f2/rel; 
#This is $\rho_4(B_3)$ with the presentation given in (A.1)
for h in Filtered( AllSubgroups(rho4B3), g -> IsNormal(rho4B3, g)) 
do Print( StructureDescription(h), "\n"); od;
    \end{lstlisting}
    and we get the following information from GAP:
    \begin{lstlisting}[language=GAP]
1
C2
C2 x C2
C2 x C2 x C2
A4
C2 x A4
A4 : C4
    \end{lstlisting}
From this computation and the information of the previous items  of this result, mainly that $Z(\rho_4(B_3))=\Zz/2\Zz$, $\rho_4(B_3)^{(1)}=A_4$ and $\rho_4(B_3)^{(2)}=V$, we obtain item (3).
\end{enumerate}

\endproof



\bibliography{braid_groups.bib}{}

@article {ACampo:1979,
    AUTHOR = {A'Campo, N.},
     TITLE = {Tresses, monodromie et le groupe symplectique},
   JOURNAL = {Comment. Math. Helv.},
  FJOURNAL = {Commentarii Mathematici Helvetici},
    VOLUME = {54},
      YEAR = {1979},
    NUMBER = {2},
     PAGES = {318--327},
      ISSN = {0010-2571},
   MRCLASS = {14D05 (32B30)},
  MRNUMBER = {535062},
MRREVIEWER = {J. A. Morrow},
       DOI = {10.1007/BF02566275},
       URL = {https://doi.org/10.1007/BF02566275},
}

@Article{Artin:1925,
 Author = {Artin, E.},
 Title = {Theorie der {Z{\"o}pfe}.},
 FJournal = {Abhandlungen aus dem Mathematischen Seminar der Universit{\"a}t Hamburg},
 Journal = {Abh. Math. Semin. Univ. Hamb.},
 ISSN = {0025-5858},
 Volume = {4},
 Pages = {47--72},
 Year = {1925},
 Language = {German},
 DOI = {10.1007/BF02950718},
 zbMATH = {2592684},
 JFM = {51.0450.01}
}

@article {JTits:1966,
    AUTHOR = {Tits, J.},
    TITLE= {Normalisateurs de tores I: Groupes de Coxeter étendus},
    FJournal = {Journal of Algebra},
 Journal = {J. Algebra},
    VOLUME={4},
    YEAR={1966}, 
    PAGES={96--116},
}

@article {Arnold:1968,
    AUTHOR = {Arnol'd, V. I.},
     TITLE = {A remark on the branching of hyperelliptic integrals as
              functions of the parameters},
   JOURNAL = {Funkcional. Anal. i Prilo\v{z}en.},
  FJOURNAL = {Akademija Nauk SSSR. Funkcional\cprime nyi Analiz i ego Prilo\v{z}enija},
    VOLUME = {2},
      YEAR = {1968},
    NUMBER = {3},
     PAGES = {1--3},
      ISSN = {0374-1990},
   MRCLASS = {14.55},
  MRNUMBER = {0244266},
MRREVIEWER = {Wazir Hasan Abdi},
}

@Article{Au,
 Author = {Auslander, L.},
 Title = {Bieberbach's theorems on space groups and discrete uniform subgroups of {Lie} groups},
 FJournal = {Annals of Mathematics. Second Series},
 Journal = {Ann. Math. (2)},
 ISSN = {0003-486X},
 Volume = {71},
 Pages = {579--590},
 Year = {1960},
 Language = {English},
 DOI = {10.2307/1969945},
 zbMATH = {3162064},
 Zbl = {0099.25602}
}

@article{Brendle-Margalit:2018,
    AUTHOR = {Brendle, T. E. and Margalit, D.},
    TITLE = {The level four braid group},
    JOURNAL = {J. Reine Angew. Math.},
    FJOURNAL = {Journal f\"{u}r die Reine und Angewandte Mathematik. [Crelle's Journal]},
    VOLUME = {735},
    YEAR = {2018},
    PAGES = {249--264},
    ISSN = {0075-4102},
    MRCLASS = {57M07 (20F36)},
    MRNUMBER = {3757477},
    MRREVIEWER = {Matthew C. B. Zaremsky},
    DOI = {10.1515/crelle-2015-0032},
    URL = {https://doi.org/10.1515/crelle-2015-0032},
}

@incollection {Gambaudo-Ghys:2005,
  AUTHOR = {Gambaudo, J.-M. and Ghys, \'{E}.},
  TITLE = {Braids and signatures},
  BOOKTITLE = {Six papers on signatures, braids and {S}eifert surfaces},
  SERIES = {Ensaios Mat.},
  VOLUME = {30},
  PAGES = {174--216},
  NOTE = {Reprinted from Bull. Soc. Math. France {{\textbf{1}}33} (2005), no. 4, 541--579 [ MR2233695]},
  PUBLISHER = {Soc. Brasil. Mat., Rio de Janeiro},
  YEAR = {2016},
  MRCLASS = {57M25 (20F36)},
  MRNUMBER = {3617348},
}

@Article{Goncalves-Guaschi:2004,
 Author = {Gon{\c{c}}alves, D.~L. and Guaschi, J.},
 Title = {The roots of the full twist for surface braid groups.},
 FJournal = {Mathematical Proceedings of the Cambridge Philosophical Society},
 Journal = {Math. Proc. Camb. Philos. Soc.},
 ISSN = {0305-0041},
 Volume = {137},
 Number = {2},
 Pages = {307--320},
 Year = {2004},
 Language = {English},
 DOI = {10.1017/S0305004104007595},
 zbMATH = {2113297},
 Zbl = {1089.20022}
}

@article {Goncalves-Guaschi-Ocampo:2017,
    AUTHOR = {Gon\c{c}alves, D.~L. and Guaschi, J. and Ocampo, O.},
     TITLE = {A quotient of the {A}rtin braid groups related to
              crystallographic groups},
   JOURNAL = {J. Algebra},
  FJOURNAL = {Journal of Algebra},
    VOLUME = {474},
      YEAR = {2017},
     PAGES = {393--423},
      ISSN = {0021-8693},
   MRCLASS = {20F36},
  MRNUMBER = {3595797},
MRREVIEWER = {Valeriy G. Bardakov},
       DOI = {10.1016/j.jalgebra.2016.11.003},
       URL = {https://doi.org/10.1016/j.jalgebra.2016.11.003},
}

@Article{BeckMarin:2020,
 Author = {Beck, V. and Marin, I.},
 Title = {Torsion subgroups of quasi-abelianized braid groups},
 FJournal = {Journal of Algebra},
 Journal = {J. Algebra},
 ISSN = {0021-8693},
 Volume = {558},
 Pages = {3--23},
 Year = {2020},
 Language = {English},
 DOI = {10.1016/j.jalgebra.2019.07.009}
}

@Article{Goncalves-Guaschi-Ocampo-Pereiro,
 Author = {Gon{\c{c}}alves, D. L. and Guaschi, J. and Ocampo, O. and Pereiro, C.},
 Title = {Crystallographic groups and flat manifolds from surface braid groups},
 FJournal = {Topology and its Applications},
 Journal = {Topology Appl.},
 ISSN = {0166-8641},
 Volume = {293},
 Pages = {16},
 Note = {Id/No 107560},
 Year = {2021},
 Language = {English},
 DOI = {10.1016/j.topol.2020.107560},
 zbMATH = {7331369},
 Zbl = {1493.20011}
}

@Article{Bellingeri-Guaschi-Makri,
 Author = {Bellingeri, P. and Guaschi, J. and Makri, S.},
 Title = {Unrestricted virtual braids and crystallographic braid groups},
 FJournal = {Bolet{\'{\i}}n de la Sociedad Matem{\'a}tica Mexicana. Third Series},
 Journal = {Bol. Soc. Mat. Mex., III. Ser.},
 ISSN = {1405-213X},
 Volume = {28},
 Number = {3},
 Pages = {16},
 Note = {Id/No 63},
 Year = {2022},
 Language = {English},
 DOI = {10.1007/s40590-022-00454-3},
 zbMATH = {7585484},
 Zbl = {1508.20039}
}

@Article{Cerqueira-Ocampo,
 Author = {Cerqueira dos Santos J{\'u}nior, P. C. and Ocampo, O.},
 Title = {Virtual braid groups, virtual twin groups and crystallographic groups},
 FJournal = {Journal of Algebra},
 Journal = {J. Algebra},
 ISSN = {0021-8693},
 Volume = {632},
 Pages = {567--601},
 Year = {2023},
 Language = {English},
 DOI = {10.1016/j.jalgebra.2023.06.005},
 zbMATH = {7710372},
 Zbl = {1521.20073}
}

@book {Dekimpe,
    AUTHOR = {Dekimpe, K.},
     TITLE = {Almost-{B}ieberbach groups: affine and polynomial structures},
    SERIES = {Lecture Notes in Mathematics},
    VOLUME = {1639},
 PUBLISHER = {Springer-Verlag, Berlin},
      YEAR = {1996},
     PAGES = {x+259},
      ISBN = {3-540-61899-6},
   MRCLASS = {20H15 (20F19 57S30)},
  MRNUMBER = {1482520},
MRREVIEWER = {Paul\ Igodt},
       DOI = {10.1007/BFb0094472},
       URL = {https://doi.org/10.1007/BFb0094472},
}

@book{Charlap,
 author = {Charlap, L. S.},
 title = {Bieberbach groups and flat manifolds},
 fseries = {Universitext},
 series = {Universitext},
 issn = {0172-5939},
 year = {1986},
 publisher = {Springer, Cham},
 language = {English},
 zbMATH = {3984018},
 Zbl = {0608.53001}
}

@book{Wolf,
 author = {Wolf, J. A.},
 title = {Spaces of constant curvature},
 edition = {6th ed.},
 isbn = {978-0-8218-5282-8},
 year = {2011},
 publisher = {Providence, RI: AMS Chelsea Publishing},
 language = {English},
 zbMATH = {5830219},
 Zbl = {1216.53003}
}

@book {Farb-Margalit:book,
    AUTHOR = {Farb, B. and Margalit, D.},
     TITLE = {A primer on mapping class groups},
    SERIES = {Princeton Mathematical Series},
    VOLUME = {49},
 PUBLISHER = {Princeton University Press, Princeton, NJ},
      YEAR = {2012},
     PAGES = {xiv+472},
      ISBN = {978-0-691-14794-9},
   MRCLASS = {57M50 (20F36 20F65 57M07 57N05)},
  MRNUMBER = {2850125},
MRREVIEWER = {Stephen P. Humphries},
}

@manual{GAP,
    key          = "GAP",
    organization = "The GAP~Group",
    title        = "{GAP -- Groups, Algorithms, and Programming,
                    Version 4.12.2}",
    year         = 2022,
    url          = "\url{https://www.gap-system.org}",
    }

@Article{GPS,
 Author = {G{\k{a}}sior, A. and Petrosyan, N. and Szczepa{\'n}ski, A.},
 Title = {Spin structures on almost-flat manifolds},
 FJournal = {Algebraic \& Geometric Topology},
 Journal = {Algebr. Geom. Topol.},
 ISSN = {1472-2747},
 Volume = {16},
 Number = {2},
 Pages = {783--796},
 Year = {2016},
 Language = {English},
 DOI = {10.2140/agt.2016.16.783},
 zbMATH = {6577066},
 Zbl = {1338.53073}
}

@Book{Johnson,
 Author = {Johnson, D. L.},
 Title = {Presentations of groups},
 FSeries = {London Mathematical Society Student Texts},
 Series = {Lond. Math. Soc. Stud. Texts},
 ISSN = {0963-1631},
 Volume = {15},
 ISBN = {0-521-37824-9; 0-521-37203-8},
 Year = {1990},
 Publisher = {Cambridge etc.: Cambridge University Press},
 Language = {English},
 zbMATH = {194089},
 Zbl = {0696.20027}
}

@book {Magnus-Karrass-Solitar:book,
    AUTHOR = {Magnus, W. and Karrass, A. and Solitar, D.},
     TITLE = {Combinatorial group theory},
   EDITION = {second},
      NOTE = {Presentations of groups in terms of generators and relations},
 PUBLISHER = {Dover Publications, Inc., Mineola, NY},
      YEAR = {2004},
     PAGES = {xii+444},
      ISBN = {0-486-43830-9},
   MRCLASS = {20E05 (05A15 20E06 20F12)},
  MRNUMBER = {2109550},
}

@Article{OR,
 Author = {Ocampo, O. and Rodr{\'{\i}}guez-Nieto, J. G.},
 Title = {On {Bieberbach} subgroups of {{\(B _n/[P _n,P _n]\)}} and flat manifolds with cyclic holonomy {{\(\mathbb{Z}_{2^d}\)}}},
 FJournal = {Topology and its Applications},
 Journal = {Topology Appl.},
 ISSN = {0166-8641},
 Volume = {265},
 Pages = {12},
 Note = {Id/No 106827},
 Year = {2019},
 Language = {English},
 DOI = {10.1016/j.topol.2019.106827},
 zbMATH = {7103808},
 Zbl = {1515.20173}
}

@Article{Porteous,
 Author = {Porteous, H. L.},
 Title = {Anosov diffeomorphisms of flat manifolds},
 FJournal = {Topology},
 Journal = {Topology},
 ISSN = {0040-9383},
 Volume = {11},
 Pages = {307--315},
 Year = {1972},
 Language = {English},
 DOI = {10.1016/0040-9383(72)90016-X},
 zbMATH = {3376478},
 Zbl = {0237.58015}
}

@article {Stylianakis:2018,
    AUTHOR = {Stylianakis, C.},
     TITLE = {Congruence subgroups of braid groups},
   JOURNAL = {Internat. J. Algebra Comput.},
  FJOURNAL = {International Journal of Algebra and Computation},
    VOLUME = {28},
      YEAR = {2018},
    NUMBER = {2},
     PAGES = {345--364},
      ISSN = {0218-1967},
   MRCLASS = {20F36 (20F05 20F65)},
  MRNUMBER = {3786423},
MRREVIEWER = {Matthew C. B. Zaremsky},
       DOI = {10.1142/S0218196718500169},
       URL = {https://doi.org/10.1142/S0218196718500169},
}

@Misc{Coxeter:1957,
 Author = {Coxeter, H. S. M.},
 Title = {Factor groups of the braid group},
 Year = {1959},
 Language = {English},
 HowPublished = {Proc. 4th {Can}. {Math}. {Congr}., {Banff} 1957, 95--122 (1959).},
 zbMATH = {3152402},
 Zbl = {0093.25003}
}

@article {Frasch:1933,
    AUTHOR = {Frasch, H.},
     TITLE = {Die Erzeugenden der Hauptkongruenzgruppen für Primzahlstufen},
   JOURNAL = {Math. Ann.},
  FJOURNAL = {Mathematische Annalen},
    VOLUME = {108},
      YEAR = {1933},
    NUMBER = {1},
     PAGES = {229--252},
      ISSN = {1432-1807},
       DOI = {10.1007/BF01452835},
       URL = {https://doi.org/10.1007/BF01452835},
}

@article{bloomquist2023quotients,
 author = {Bloomquist, W. and Patzt, P. and Scherich, N.},
 title = {Quotients of braid groups by their congruence subgroups},
 fjournal = {Proceedings of the American Mathematical Society. Series B},
 journal = {Proc. Am. Math. Soc., Ser. B},
 issn = {2330-1511},
 volume = {11},
 pages = {508--524},
 year = {2024},
 language = {English},
 doi = {10.1090/bproc/200},
 zbMATH = {7939033},
 Zbl = {1552.20169}
}

@article{ABGH,
  title={On quotients of congruence subgroups of braid groups.},
  author={Appel, J. and Bloomquist, W. and Gravel, K. and Holden, A.},
  journal={arXiv: Group Theory},
  year={2020}
}

@incollection {BirmanBrendle:2005,
    AUTHOR = {Birman, J. S. and Brendle, T. E.},
     TITLE = {Braids: a survey},
 BOOKTITLE = {Handbook of knot theory},
     PAGES = {19--103},
 PUBLISHER = {Elsevier B. V., Amsterdam},
      YEAR = {2005},
   MRCLASS = {57M25 (20F36 55R80 57M27 57R17)},
  MRNUMBER = {2179260 (2007a:57004)},
MRREVIEWER = {Darren D. Long},
       DOI = {10.1016/B978-044451452-3/50003-4},
       URL = {http://dx.doi.org/10.1016/B978-044451452-3/50003-4},
}

@article{BDOS:2024,
 author = {Bellingeri, P. and Damiani, C. and Ocampo, O. and Stylianakis, C.},
 title = {Congruence subgroups of braid groups and crystallographic quotients. {I}},
 fjournal = {Journal of the Australian Mathematical Society},
 journal = {J. Aust. Math. Soc.},
 issn = {1446-7887},
 volume = {118},
 number = {2},
 pages = {145--163},
 year = {2025},
 language = {English},
 doi = {10.1017/S1446788724000089},
 zbMATH = {8018268}
}

@article{Tuba-Wenzl:2000,
author = {Tuba, I. and Wenzl, H.},
Title = {Representations of the braid group {{\(B_3\)}} and of {{\(\text{SL}(2,\mathbb{Z})\)}}.},
 FJournal = {Pacific Journal of Mathematics},
 Journal = {Pac. J. Math.},
 ISSN = {1945-5844},
 Volume = {197},
 Number = {2},
 Pages = {491--510},
 Year = {2001},
 Language = {English},
 DOI = {10.2140/pjm.2001.197.491},
 zbMATH = {1621006},
 Zbl = {1056.20025}
}

@article{LIMAGONCALVES2019160,
title = {Almost-crystallographic groups as quotients of Artin braid groups},
FJournal = {Journal of Algebra},
 Journal = {J. Algebra},
volume = {524},
pages = {160-186},
year = {2019},
issn = {0021-8693},
doi = {https://doi.org/10.1016/j.jalgebra.2019.01.010},
url = {https://www.sciencedirect.com/science/article/pii/S0021869319300341},
author = {Gon\c{c}alves, D.~L. and Guaschi, J. and Ocampo, O.}
}

@Misc{Banerjee-Huxford,
 Author = {Banerjee, I. and Huxford, P.},
 Title = {Generators for the level $m$ congruence subgroups of braid groups},
 Year = {2024},
 HowPublished = {Preprint, {arXiv}:2409.09612 [math.{GR}] (2024)},
 URL = {https://arxiv.org/abs/2409.09612},
 arXiv = {arXiv:2409.09612}
}

@article{Nakamura:2021,
  title={The cohomology class of the mod 4 braid group},
  author={T. Nakamura},
  JOURNAL = {Internat. J. Algebra Comput.},
  FJOURNAL = {International Journal of Algebra and Computation},
  year={2021},
  volume={32},
  pages={1125-1159},
  url={https://api.semanticscholar.org/CorpusID:236428456}
}

@article{KordekMargalit:2022,
  title={Representation stability in the level 4 braid group},
  author={Kordek, K. and Margalit, D.},
   FJournal = {Mathematische Zeitschrift},
 Journal = {Math. Z.},
  year={2022},
  volume={302},
  pages={433–482},
  doi={https://doi.org/10.1007/s00209-022-03059-8},
  issn = {1432-1823} 
}

@article{Bruillard-Plavnik-Rowell:2015,
    author = {Bruillard, P. and Chang, L. and Hong, S. M. and Plavnik, J. Y. and Rowell, E. C. and Sun, M. Y.},
 Title = {Low-dimensional representations of the three component loop braid group.},
 FJournal = {Journal of Mathematical Physics},
 Journal = {J. Math. Phys.},
 ISSN = {0022-2488},
 Volume = {56},
 Number = {11},
 Pages = {111707, 15},
 Year = {2015},
 Language = {English},
 DOI = {10.1063/1.4935361},
 URL = {hdl.handle.net/11336/51849},
 zbMATH = {6519491},
 Zbl = {1331.20045}
}

@book {Fulton-Harris:1991,
    AUTHOR = {Fulton, W. and Harris, J.},
     TITLE = {Representation theory},
    SERIES = {Graduate Texts in Mathematics},
    VOLUME = {129},
      NOTE = {A first course,
              Readings in Mathematics},
 PUBLISHER = {Springer-Verlag, New York},
      YEAR = {1991},
     PAGES = {xvi+551},
      ISBN = {0-387-97527-6; 0-387-97495-4},
   MRCLASS = {20G05 (17B10 20G20 22E46)},
  MRNUMBER = {1153249},
MRREVIEWER = {James\ E.\ Humphreys},
       DOI = {10.1007/978-1-4612-0979-9},
       URL = {https://doi.org/10.1007/978-1-4612-0979-9},
}
\bibliographystyle{alpha}

\end{document}